\documentclass[12pt,final]{article}
\usepackage{amsmath,amssymb,amsfonts}
\usepackage[utf8]{inputenc} 
\usepackage[T1]{fontenc}
\usepackage{enumerate}
\usepackage{pdfpages}
\usepackage{hyperref}
\usepackage{graphicx}
\usepackage{tabularx}
\usepackage{url}
\newcommand{\MyParagraph}[1]{\bigskip \noindent \textbf{#1.}\;}
\newcommand{\MyCite}[1]{\colred{References}}
\newtheorem{Theorem}{Theorem}
\newcommand{\NN}{\mathbb{N}}
\newcommand{\ZZ}{\mathbb{Z}}
\newcommand{\RR}{\mathbb{R}}
\newcommand{\CC}{\mathbb{C}}
\newcommand{\dd}{\ensuremath{\mathrm{d}}}
\newcommand{\ee}{\ensuremath{\mathrm{e}}}
\newcommand{\ii}{\ensuremath{\mathrm{i}}}
\newcommand{\nE}{\mathcal{E}}
\newcommand{\nF}{\mathcal{F}}
\newcommand{\nS}{\mathcal{S}}
\newcommand{\nL}{\mathcal{L}}
\newcommand{\nO}{\mathcal{O}}
\newcommand{\nSup}[1]{\nS^{(#1)}}
\newcommand{\abs}[1]{|#1|}
\newcommand{\absbig}[1]{\big|#1\big|}
\newcommand{\absBig}[1]{\Big|#1\Big|}
\newcommand{\norm}[2]{\|#2\|_{#1}}
\newcommand{\normbig}[2]{\big\|#2\big\|_{#1}}
\newcommand{\scp}[3]{(#2\vert#3)_{#1}}
\newcommand{\taun}{\tau_{n}}
\begin{document}
\title
{\textbf{Splitting methods with complex coefficients for linear and nonlinear evolution equations}}
\author
{S.~Blanes, F.~Casas, C.~Gonz{\'a}lez, M.~Thalhammer%
\footnote{Addresses: Sergio Blanes, Universitat Polit{\`e}cnica de Val{\`e}ncia, Instituto de Matem{\'a}tica Multidisciplinar, 46022~Valencia, Spain. 
Fernando Casas, Universitat Jaume~I, IMAC and Departament de Matem{\`a}tiques, 12071~Castell\'on, Spain. 
Ces{\'a}reo Gonz{\'a}lez, Universidad de Valladolid, Departamento de Matem{\'a}tica Aplicada, 47011~Valladolid, Spain. 
Mechthild Thalhammer, Universit{\"a}t Innsbruck, Institut f{\"u}r Mathematik, 6020~Innsbruck, Austria.
Email addresses: \url{serblaza@imm.upv.es}, \url{fernando.casas@uji.es}, \url{ome@am.uva.es}, \url{mechthild.thalhammer@uibk.ac.at}.
Websites: \url{personales.upv.es/serblaza}, \url{www.gicas.uji.es/fernando.html}, \url{www.imuva.uva.es/en/investigadores/33}, \url{techmath.uibk.ac.at/mecht}}
}
\date{}
\maketitle
\begin{abstract}
\noindent 
This contribution is dedicated to the exploration of exponential operator splitting methods for the time integration of evolution equations.
It entails the review of previous achievements as well as the depiction of novel results.
The standard class of splitting methods involving real coefficients is contrasted with an alternative approach that relies on the incorporation of complex coefficients. 
In view of long-term computations for linear evolution equations, it is expedient to distinguish symmetric, symmetric-conjugate, and alternating-conjugate schemes.
The scope of applications comprises high-order reaction-diffusion equations and complex Ginzburg--Landau equations, which are of relevance in the theories of patterns and superconductivity. 
Time-dependent Gross--Pitaevskii equations and their parabolic counterparts, which model the dynamics of Bose–Einstein condensates and arise in ground state computations, are formally included as special cases. 
Numerical experiments confirm the validity of theoretical stability conditions and global error bounds as well as the benefits of higher-order complex splitting methods in comparison with standard schemes.
\end{abstract}
\section{Introduction}
\MyParagraph{Splitting methods}
Exponential operator splitting methods constitute nowadays a popular class of numerical algorithms for the integration in time of initial value problems
\begin{subequations}
\label{eq:Introduction}
\begin{equation}
\begin{cases}
u'(t) = F\big(u(t)\big)\,, \quad t \in [t_0, T]\,, \\
u(t_0) \text{ given}\,, 
\end{cases}
\end{equation}
involving ordinary or partial differential equations, respectively. 
A fundamental premise is that the defining function naturally decomposes into two or more parts 
\begin{equation}
F = \sum_{\ell = 1}^{L} F_{\ell}\,, 
\end{equation}
\end{subequations}
so that the problem under consideration can be subdivided into different initial value problems that are easier to solve than the original problem.
For detailed information on splitting methods, we refer to~\cite{HairerLubichWanner2006,McLachlanQuispel2002,SanzSernaCalvo2018} and the recent review~\cite{BlanesCasasMurua2024}.

\MyParagraph{Benefits}
In many relevant instances, exponential operator splitting methods are straightforward to implement and present favourable properties concerning error propagation.
In addition, they preserve a variety of structural properties possessed by the differential equation.
As a consequence, splitting methods are widely used in computational quantum mechanics when dealing with problems requiring to solve time-dependent Schr{\"o}dinger equations, see~\cite{Gross1961,Messiah1999,Pitaevskii1961,PitaevskiiStringari2003} as well as~\cite{Bao2004,BaoJinMarkowich2002,BesseBidegarayDescombes2002,BlanesCasasMurua2015,CaliariZuccher2021,CarlesSu2024,IserlesKropielnicka2024,JahnkeLubich2000,
Thalhammer2008,Thalhammer2012,ThalhammerAbhau2012} and references given therein. 

\MyParagraph{Limitations}
However, the applicability of exponential operator splitting methods to non-reversible systems has been more limited, see for example~\cite{BertoliVilmart2020,TapleyCelledoniOwrenAndersson2019}.
One important reason for that is related with an essential feature of this class of time integration methods: splitting methods necessarily involve negative coefficients when their order is three or higher, see~\cite{AuzingerHofstaetterKoch2019,BlanesCasas2005}.
As a consequence, severe instabilities arise when higher-order splitting methods are applied to evolution equations of parabolic type such as reaction-diffusion equations.

\MyParagraph{Alternative approach}
An alternative to deal with higher-order exponential operator splitting methods for equations evolved by semigroups was proposed in \cite{BandraukShen1991,CastellaChartierDecombesVilmart2009,HansenOstermann2009}, see also~\cite{EngelNagel2000,HillePhillips1974,Lunardi2013,Pazy1983}: using instead \emph{complex} coefficients having \emph{positive real part}.
In that way, the stability is recovered  in certain situations, amongst others, for reaction-diffusion systems involving real constants.
Moreover, increasing the order of the scheme typically improves its efficiency over a wide range of accuracies, in contrast with what takes place with real coefficients, see~\cite{BlanesCasasChartierMurua2013}.
In view of these appealing properties, the practical implementation and systematic exploration of exponential operator splitting methods with complex coefficients has been carried out, not only by designing new high-order schemes but also by analysing their preservation properties, for both ordinary and partial differential equations \cite{BaderBlanesCasas2013,BernierBlanesCasasEscorihuela2023,BlanesCasasChartierEscorihuela2022,BlanesCasasEscorihuela2022,BlanesCasasGonzalezThalhammer2024}.

\MyParagraph{Additional features}
When considering standard exponential operator splitting methods with real coefficients, left-right palindromic or so-called symmetric compositions are usually preferred, since their construction is easier and in addition they provide time-symmetric approximations to the exact solution.
It turns out, however, that in the complex case other alternatives are more favourable with respect to the preservation of properties.
The benefits of so-called symmetric-conjugate schemes which are symmetric in the real part and antisymmetric in the imaginary part have been illustrated in~\cite{BernierBlanesCasasEscorihuela2023,BlanesCasasChartierEscorihuela2022} for linear unitary problems and in~\cite{BlanesCasasGonzalezThalhammer2024} for parabolic equations. 
Furthermore, it is promising to study complex splitting methods that are obtained by concatenating a given composition with its complex conjugate.
Work in progress~\cite{BernierBlanesCasasEscorihuela2024} indicates that schemes of this kind exhibit an excellent behaviour in long-term computations of linear evolution equations.

\MyParagraph{Objective}
In this contribution, our main objective is to analyse the performance of high-order exponential operator splitting methods involving complex coefficients for a general model problem that is defined by a linear combination of powers of the Laplace operator, a space-dependent function, and a nonlinear multiplication operator.
We deduce stability conditions and state a convergence result which generalises our former analysis of standard real splitting methods within the context of linear and nonlinear Schr{\"o}dinger equations, see~\cite{Thalhammer2008,Thalhammer2012}.
The scope of applications includes reaction-diffusion and Ginzburg--Landau-type equations, which are of relevance, amongst others, in the theories of patterns and superconductivity.  
More precisely, specifically designed high-order reaction-diffusion equations reveal the formation of quasicrystalline patterns, and complex Ginzburg--Landau equations describe the dynamical behaviour of generic spatially extended systems undergoing a supercritical Hopf bifurcation from a stationary to an oscillatory state, see~\cite{JiangZhang2014,Saarloos1995}.
Time-dependent Gross--Pitaevskii equations modelling the dynamics of Bose–Einstein condensates and their parabolic counterparts arising in ground and excited state computations are formally retained as particular cases, see~\cite{Bao2004,Gross1961,Pitaevskii1961,PitaevskiiStringari2003}.

\MyParagraph{Outline}
This contribution has the following structure. 
In Section~\ref{sec:SplittingMethods}, we introduce a general framework of abstract evolution equations and specify the format of complex exponential operator splitting methods.
In Section~\ref{eq:ModelProblems}, we detail our model problem including high-order reaction-diffusion equations, complex Ginzburg--Landau equations, time-dependent Gross--Pitaevskii equations, and their parabolic counterparts as particular instances.  
Important means for a stability and error analysis of high-order exponential operator splitting methods are summarised in Section~\ref{sec:Convergence}. 
Hereby, it is practical to include instructions along the lines of~\cite{Thalhammer2008}, where standard real splitting methods have been studied for linear Schr{\"o}dinger equations comprising the Laplace operator and a regular space-dependent function. 
In essence, the extension to the significantly more involved nonlinear case relies on the powerful calculus of Lie derivatives, see~\cite{Varadarajan1984}.
The rigorous treatment of full discretisations combining time-splitting methods and spectral space discretisations requires suitable adaptations of~\cite{Thalhammer2012}.
Numerical experiments that confirm the theoretical statements and the advantages of high-order complex splitting methods in comparison with standard schemes are finally presented in Section~\ref{sec:NumericalExperiments}.
\section{Splitting methods}
\label{sec:SplittingMethods}
In this section, we collect basic information on exponential operator splitting methods. 
For a thorough description of the analytical frameworks and a thematic classification of splitting methods within the field of geometric integration, we refer to~\cite{EngelNagel2000,HillePhillips1974,Lunardi2013,Pazy1983} and~\cite{BlanesCasasMurua2024,HairerLubichWanner2006,McLachlanQuispel2002,SanzSernaCalvo2018}.
  
\MyParagraph{Abstract evolution equations}
The starting point of our considerations are initial value problems for abstract evolution equations that can be cast into the general format 
\begin{equation}
\label{eq:IVP}
\begin{cases}
u'(t) = F\big(u(t)\big) = F_1\big(u(t)\big) + F_2\big(u(t)\big)\,, \quad t \in [t_0, T]\,, \\
u(t_0) \text{ given}\,,
\end{cases}
\end{equation}
see also~\eqref{eq:Introduction}.
Throughout, we denote by $(X, \norm{X}{\cdot})$ the underlying Banach spaces with associated norms and assume that the domains of the defining operators $F_{\ell}: D(F_{\ell}) \subseteq X \to X$, $\ell \in \{1, 2\}$, have non-empty intersections.
Within the context of nonlinear evolution equations of parabolic or Schr{\"o}dinger type, it is expedient to employ the frameworks of analytical semigroups or strongly continuous unitary groups, respectively.

\MyParagraph{Splitting approach}
Exponential operator splitting methods for an evolution equation of the form~\eqref{eq:IVP} rely on the presumption that the numerical approximation of the associated subproblems 
\begin{equation}
\label{eq:IVPSubproblems}
\begin{cases}
u_j'(t) = F_j\big(u_j(t)\big)\,, \quad t \in [t_0, T]\,, \\
u_j(t_0) \text{ given}\,, 
\end{cases} \quad j \in \{1,2\}\,,
\end{equation}
is significantly simpler compared to the numerical approximation of the original problem.

\MyParagraph{Scope of applications}
The scope of applications includes Hamiltonian systems from classical mechanics and Schr{\"o}dinger equations from quantum mechanics, where the advantages of exponential operator splitting methods constituting geometric numerical integrators become apparent.
It naturally extends to non-reversible systems such as reaction-diffusion systems and Ginzburg--Landau-type equations, which form beautiful spatio-temporal patterns and attract interest, amongst others, in biology, chemistry, geology, and physics.
Moreover, the fundamental concepts are applicable to higher-order damped wave equations from nonlinear acoustics and kinetic equations from plasma physics.
See for instance~\cite{BlanesCasasMurua2024,CasasCrouseillesFaouMehrenberger2017,JiangZhang2014,KaltenbacherNikolicThalhammer2015,LiEtAl2019,Saarloos1995} and references given therein.

\MyParagraph{Model problems}
We study a concrete model problem that comprises a linear combination of powers of the Laplace operator, a space-dependent function, and a nonlinear multiplication operator. 
Even though the qualitative characteristics of solutions are different, high-order reaction-diffusion equations describing quasicrystals as well as complex Ginzburg--Landau equations arising in superconductivity are formally included.
Time-dependent Gross--Pitaevskii equations modelling the dynamics of Bose--Einstein-condensates and their parabolic counterparts emerging in ground and excited state computations by imaginary time propagation are retained as special cases.
The reductions to related linear cases lead to partial differential equations that involve, on the one hand, linear differential operators and, on the other hand, multiplication operators defined by space-dependent functions. 

\MyParagraph{Evolution operators}
Henceforth, we use a symbolic notation for the exact evolution operators associated with the original problem~\eqref{eq:IVP} or the subproblems~\eqref{eq:IVPSubproblems}, respectively, and indicate the dependencies on the current times, the defining operators, and the initial states.
For instance, we set 
\begin{equation*}
\nE_{t, F}\big(u(t_0)\big) = u(t)\,, \quad t \in [t_0, T]\,,
\end{equation*}
for the evolution operator corresponding to~\eqref{eq:IVP}.

\MyParagraph{Iterated commutators}
Under the requirement of suitably restricted domains, the commutator of two nonlinear operators is defined by
\begin{subequations}
\label{eq:IteratedCommutators}
\begin{equation}
\big[F_1, F_2\big](v) = F_1'(v) \, F_2(v) - F_2'(v) \, F_1(v)\,.
\end{equation}  
More generally, the iterated commutators are determined recursively
\begin{equation}
\text{ad}_{F_1}^{\, \ell}(F_2) = \begin{cases} F_2\,, & \ell = 0\,, \\ \big[F_1, \text{ad}_{F_1}^{\, \ell - 1}(F_2)\big]\,, &\ell \in \NN_{\geq 1}\,. \end{cases}
\end{equation}
\end{subequations}

\MyParagraph{Complex splitting methods}
We employ a usual time-stepping approach, denoting by  
\begin{subequations}
\label{eq:Splitting}  
\begin{equation}
t_0 < t_1 < \dots < t_N = T\,, \quad \taun = t_{n+1} - t_n\,, \quad n \in \{0, 1, \dots, N-1\}\,,
\end{equation}
the underlying time grid points with corresponding time stepsizes.
The recurrence relations for the exact and numerical solution values read as 
\begin{equation}
\begin{gathered}
u_{n+1} = \nS_{\taun,F}(u_n) \approx u(t_{n+1}) = \nE_{\taun, F}\big(u(t_n)\big)\,, \\
n \in \{0, 1, \dots, N-1\}\,. 
\end{gathered}
\end{equation}
For suitable choices of the defining method coefficients, complex exponential operator splitting methods can be written as  
\begin{equation}
\label{eq:SplittingS}  
\begin{gathered}
\nS_{\taun,F} = \nE_{\taun, b_s F_2} \circ \, \nE_{\taun, a_s F_1} \circ \dots \circ \, \nE_{\taun, b_1 F_2} \circ \, \nE_{\taun, a_1 F_1}\,, \\
(a_j, b_j)_{j=1}^{s} \in \CC^{2s}\,.
\end{gathered}
\end{equation}
\end{subequations}
Evidently, standard splitting methods are included as special cases with real coefficients.

\MyParagraph{Order and error estimates}
Within the context of nonstiff differential equations with sufficiently regular solutions, the order $p \in \NN_{\geq 1}$ of a complex exponential operator splitting method~\eqref{eq:Splitting} is defined through the relation
\begin{equation}
\label{eq:LocalErrorOrder}
\nL_{\taun,F}(v) = \nS_{\taun,F}(v) - \nE_{\taun, F}(v) = \nO\big(\taun^{p+1}\big) \quad \text{as} \quad \taun \searrow 0\,.
\end{equation}

\MyParagraph{Global error estimates}
In view of a rigorous convergence analysis of exponential operator splitting methods, we employ a standard argument based on the telescopic identity to conclude that the validity of stability bounds combined with suitable local error expansions ensuring~\eqref{eq:LocalErrorOrder} implies global error estimates 
\begin{equation*}
\begin{gathered}  
\normbig{X}{u_n - u(t_n)} \leq C \, \Big(\normbig{X}{u_0 - u(t_0)} + \tau_{\max}^p\Big)\,, \\
n \in \{1, \dots, N\}\,, \quad \tau_{\max} = \max_{n \in \{0, 1, \dots, N-1\}} \taun\,. 
\end{gathered}
\end{equation*}
In order to subsume our modus operandi, we meanwhile restrict ourselves to evolution equations~\eqref{eq:IVP} that are defined by bounded linear operators 
\begin{equation*}
F_1, F_2: X \longrightarrow X 
\end{equation*}
and complex splitting methods~\eqref{eq:Splitting} applied with constant time stepsizes
\begin{equation*}
\taun = \tau\,, \quad n \in \{0, 1, \dots, N-1\}\,, 
\end{equation*}
since then the exact and numerical solutions can be represented by exponential series
\begin{equation*}  
\nE_{\tau, F} = \ee^{\, \tau \, (F_1 + F_2)}\,, \quad \nS_{\tau, F} = \ee^{\, b_s \tau F_2} \, \ee^{\, a_s \tau F_1} \cdots \, \ee^{\, b_1 \tau F_2} \, \ee^{\, a_1 \tau F_1}\,.
\end{equation*}
Consequently, stability bounds follow at once from 
\begin{equation*}  
\begin{gathered}
\normbig{X \leftarrow X}{\nE_{\tau, F}} \leq \ee^{\, C \tau}\,, \quad \normbig{X \leftarrow X}{\nS_{\tau, F}} \leq \ee^{\, C \tau}\,, \\
C = \max\Bigg\{\norm{X \leftarrow X}{F_1} + \norm{X \leftarrow X}{F_2}, \sum_{j=1}^{s} \Big(\abs{a_j} \norm{X \leftarrow X}{F_1} + \abs{b_j} \norm{X \leftarrow X}{F_2}\Big)\bigg\}\,.
\end{gathered}
\end{equation*}
Provided that the method coefficients satisfy certain order conditions such that local error expansions~\eqref{eq:LocalErrorOrder} hold true, the desired global error estimates are obtained by means of the telescopic identity 
\begin{equation*}
\begin{gathered}
u_n - u(t_n) = \big(\nS_{\tau, F}\big)^n \, \big(u_0 - u(t_0)\big) 
+ \sum_{\ell=0}^{n-1} \big(\nS_{\tau, F}\big)^{n-1-\ell} \nL_{\tau, F} \, \big(\nE_{\tau, F}\big)^{\ell} \, u(t_0)\,, \\
n \in \{1, \dots, N\}\,, 
\end{gathered}
\end{equation*}
and straightforward estimation
\begin{equation*}
\begin{gathered}
\normbig{X}{u_n - u(t_n)} \leq \ee^{\, C t_n} \, \Big(\normbig{X}{u_0 - u(t_0)} + n \, \normbig{X \leftarrow X}{\nL_{\tau,F}} \normbig{X}{u(t_0)}\Big)\,, \\
n \in \{1, \dots, N\}\,. 
\end{gathered}
\end{equation*}
Theorem~\ref{thm:Theorem} given in Section~\ref{sec:Convergence} below constitutes the natural generalisation to abstract evolution equations involving unbounded nonlinear operators. 

\MyParagraph{Additional features}
Whenever the primary interest of the investigations concerns the specification of exponential operator splitting methods, it is convenient to supress the dependencies on the time stepsizes as well as the defining operators and to replace~\eqref{eq:SplittingS} by the compact representation
\begin{equation}
\label{eq:SShort}
\nS = \big(b_s, a_s, \dots, b_1, a_1\big)\,.
\end{equation}
As special patterns in the distribution of the method coefficients are of importance within the context of long-term computations, we distinguish symmetric (s), symmetric-conjugate (sc), and alternating-conjugate (ac) schemes.
Assuming for a moment that the defining operators correspond to real symmetric matrices $F_1, F_2 \in \RR^{M \times M}$, these properties can be connected to complex conjugation and transposition of certain components, e.g.  
\begin{subequations}
\label{eq:Structures}  
\begin{equation}
\begin{gathered}
\nS_0 = \big(b_r, a_r, \dots, b_1, a_1\big)\,, \\
\overline{\nS_0} = \big(\overline{b_r}, \overline{a_r}, \dots, \overline{b_1}, \overline{a_1}\big)\,, \quad
\nS_0^T = \big(a_1, b_1, \dots, a_r, b_r\big)\,.
\end{gathered}
\end{equation}
That is, when setting $a_1 = 0$, we impose the following configurations 
\begin{equation}
\begin{gathered}
\nS_0 = \big(b_r, a_r, \dots, b_2, a_2, b_1\big)\,, \\
\nSup{\text{s}} = \big(\nS_0^T, a_{r+1}, \nS_0, 0\big)\,, \\
a_{r+1} \in \RR\,, \quad \nSup{\text{sc}} = \big(\overline{\nS_0}^T, a_{r+1}, \nS_0, 0\big)\,, \\
a_{r+1} \in \RR\,, \quad \nS_1 = \big(\overline{\nS_0}^T, a_{r+1}, \nS_0\big)\,, \quad \nSup{\text{ac}} = \big(\overline{\nS_1}, \nS_1, 0\big)\,.
\end{gathered}
\end{equation}
For $b_s = 0$, we instead obtain 
\begin{equation}
\begin{gathered}
\nS_0 = \big(a_r, \dots, b_2, a_2, b_1, a_1\big)\,, \\
\nSup{\text{s}} = \big(0, \nS_0^T, b_r, \nS_0\big)\,, \\
b_r \in \RR\,, \quad \nSup{\text{sc}} = \big(0, \overline{\nS_0}^T, b_r, \nS_0\big)\,, \\
b_r \in \RR\,, \quad \nS_1 = \big(\overline{\nS_0}^T, b_r, \nS_0\big)\,, \quad \nSup{\text{ac}} = \big(0, \overline{\nS_1}, \nS_1\big)\,.
\end{gathered}
\end{equation}
\end{subequations}
Further information is found in~\cite{BernierBlanesCasasEscorihuela2023,BernierBlanesCasasEscorihuela2024,BlanesCasasChartierEscorihuela2022,BlanesCasasGonzalezThalhammer2024}. 

\MyParagraph{Examples}
Widely used low-order exponential operator splitting methods such as the non-symmetric first-order Lie--Trotter splitting method and the symmetric second-order Strang splitting method are retained for the particular choices 
\begin{equation*}
\begin{gathered}  
\text{Lie--Trotter}: \quad p = s = 1\,, \quad \nS = (1, 1)\,, \\
\text{Strang}: \quad p = s = 2\,, \quad \nS = \big(\tfrac{1}{2}, 1, \tfrac{1}{2}, 0\big)\,,
\end{gathered}  
\end{equation*}
see~\eqref{eq:Splitting}, \eqref{eq:LocalErrorOrder}, and~\eqref{eq:SShort}.
The method coefficients of various higher-order exponential operator splitting methods are collected in a \textsc{Matlab} code that is available at \url{doi.org/10.5281/zenodo.13834638}, see also Section~\ref{sec:NumericalExperiments}.
Their additional features and stability properties are displayed in Table~\ref{tab:Table1}.

\MyParagraph{Adaptivity}
For non-reversible systems such as reaction-diffusion equations and associated linear equations, see~\eqref{eq:ReactionDiffusion} below, it is desirable to enhance the efficiency of long-term simulations by incorporating a time stepsize control.
Further details are found in~\cite{BlanesCasasGonzalezThalhammer2024,BlanesCasasThalhammer2019}, where two different strategies have been exploited, see also the references given therein. 
\section{Model problems}
\label{eq:ModelProblems}
In this section, we first introduce our general model problem and then state high-order reaction-diffusion equations and complex Ginzburg--Landau-type equations that are included as special cases. 
Fundamental questions such as well-posedness are treated in~\cite{EngelNagel2000,HillePhillips1974,Lunardi2013,Pazy1983}.

\MyParagraph{General model problem}
For a space-time-dependent complex-valued solution $U: \Omega \times [t_0, T] \subset \RR^{d} \times \RR \to \CC$, we consider the model problem 
\begin{equation}
\label{eq:GeneralModel1}
\begin{cases}
\displaystyle \partial_t U(x, t) = \sum_{k=0}^{K} \alpha_k \, \Delta^k \, U(x,t) + W(x) \, U(x, t) + f\big(U(x, t)\big)\,, \\
U(x, t_0) = U_0(x)\,, \quad (x, t) \in \Omega \times [t_0, T]\,, 
\end{cases} 
\end{equation}
which comprises a linear combination of powers of the Laplace operator, a sufficiently regular space-dependent function $W: \Omega \to \CC$, and a nonlinear multiplication operator defined by the complex function $f: \CC \to \CC$.
Throughout, we assume that the complex constants $(\alpha_k)_{k=0}^K \in \CC^{K+1}$ are chosen such that well-posedness of the problem is ensured.

\MyParagraph{Compact reformulation}
In order to rewrite~\eqref{eq:GeneralModel1} in the compact form~\eqref{eq:IVP}, we set $u(t) = U(\cdot, t)$ for $t \in [t_0, T]$ and define the associated operators for any sufficiently regular space-dependent function $v: \Omega \to \CC$ through  
\begin{equation}
\label{eq:GeneralModel2}
\begin{gathered}
\big(F_1(v)\big)(x) = \sum_{k=0}^{K} \alpha_k \, \Delta^k \, v(x)\,, \\ 
\big(F_2(v)\big)(x) = W(x) \, v(x) + f\big(v(x)\big)\,, \\
x \in \Omega\,.
\end{gathered}
\end{equation}
\subsection{High-order reaction-diffusion equations}
We study high-order reaction-diffusion equations of the form 
\begin{equation}
\label{eq:ReactionDiffusion}  
\begin{cases}
\displaystyle \partial_t U(x, t) = \sum_{k=0}^{4} \alpha_k \, \Delta^k \, U(x,t) + \sum_{k=1}^{3} \beta_k \, \big(U(x, t)\big)^k\,, \\
U(x, t_0) = U_0(x)\,, \quad (x, t) \in \Omega \times [t_0, T]\,, 
\end{cases} 
\end{equation}
retained from~\eqref{eq:GeneralModel1} for $W = 0$ and a cubic polynomial 
\begin{equation*}
f: \RR \longrightarrow \RR: v \longmapsto \sum_{k=1}^{3} \beta_k \, v^k\,.
\end{equation*}
With regard to the simulation of quasicrystalline patterns in two space dimensions, we assume that the underlying spatial domains are given by Cartesian products of bounded intervals, tacitly impose periodic boundary conditions, and focus on real constants $(\alpha_k)_{k=0}^4 \in \RR^5$ and $(\beta_k)_{k=1}^3 \in \RR^3$. 
Under the basic assumption $\alpha_4 < 0$, the problem is well-posed, and the solution $U: \Omega \times [t_0, T] \to \RR$ inherits the regularity of the prescribed initial state.  
\subsection{Complex Ginzburg--Landau-type equations}
A particularity of Ginzburg--Landau-type equations
\begin{equation}
\label{eq:GinzburgLandauType}   
\begin{cases}
\partial_t U(x, t) = \alpha_1 \, \Delta U(x, t) + \alpha_0 \, U(x, t) \\
\qquad\qquad\qquad + \, \beta_1 V(x) \, U(x, t) + \beta_2 \, \absbig{U(x, t)}^{\beta_3} \, \big(U(x, t)\big)^{\beta_4}\,, \\
U(x, t_0) = U_0(x)\,, \quad (x, t) \in \Omega \times [t_0, T]\,,
\end{cases}
\end{equation}
is the interaction of the Laplace operator and the non-holomorphic nonlinearity in the presence of complex constants.
Evidently, the model problem~\eqref{eq:GinzburgLandauType}  can be cast into the general form~\eqref{eq:GeneralModel1} when setting $\alpha_k = 0$ for $j \in \{2, 3, 4\}$, $W = \beta_1 V$, and 
\begin{equation*}
f: \CC \longrightarrow \CC: v \longmapsto \beta_2 \, \abs{v}^{\beta_3} \, v^{\beta_4}\,.
\end{equation*}
In view of relevant applications, we once again suppose that the underlying spatial domains are given by Cartesian products of (sufficiently large) bounded intervals and tacitly impose periodic boundary conditions.
Likewise, we conjecture sufficiently regular solutions $U: \Omega \times [t_0, T] \to \CC$ for appropriately prescribed potentials and initial states.  
On the one hand, in order to ensure that the evolution equation in~\eqref{eq:GinzburgLandauType} is of parabolic type, we employ the restriction $\Re(\alpha_1) > 0$.
On the other hand, it is instructive to observe that nonlinear Schr{\"o}dinger equations such as Gross--Pitaevskii equations are included in~\eqref{eq:GinzburgLandauType} in the limiting case $\Re(\alpha_1) \to 0$ and $\Im(\alpha_1) \neq 0$.

\MyParagraph{Special cases}
We next specify complex Ginzburg--Landau equations, time-dependent Gross--Pitaevskii equations, and their parabolic counterparts.
In all cases, we focus on the special choices $\beta_3 = 2$ and $\beta_4 = 1$. 
\begin{enumerate}[(i)]
\item   
Complex Ginzburg--Landau equations
\begin{equation}
\label{eq:CGL}   
\begin{cases}
\partial_t U(x, t) = \alpha_1 \, \Delta U(x, t) + \alpha_0 \, U(x, t) + \beta_2 \, \absbig{U(x, t)}^2 \, U(x, t)\,, \\
U(x, t_0) = U_0(x)\,, \quad (x, t) \in \Omega \times [t_0, T]\,,
\end{cases}
\end{equation}  
are obtained from~\eqref{eq:GinzburgLandauType} for $\beta_1 = 0$.
\item 
Time-dependent Gross--Pitaevskii equations
\begin{equation}
\label{eq:GPE}   
\begin{cases}
\ii \, \partial_t U(x, t) = \alpha \, \Delta U(x, t) + \beta \, V(x) \, U(x, t) + \vartheta \, \absbig{U(x, t)}^2 \, U(x, t)\,, \\
U(x, t_0) = U_0(x)\,, \quad (x, t) \in \Omega \times [t_0, T]\,,
\end{cases}
\end{equation}  
result from~\eqref{eq:GinzburgLandauType} for the choices $\alpha_1 = - \, \ii \, \alpha$, $\alpha_0 = 0$, $\beta_1 = - \, \ii \, \beta$, and $\beta_2 = - \, \ii \, \vartheta$ with real constants $\alpha, \beta, \vartheta \in \RR$.
\item 
Related nonlinear parabolic equations 
\begin{equation}
\label{eq:GPEParabolic}   
\begin{cases}
\partial_t U(x, t) = \alpha \, \Delta U(x, t) + \beta \, V(x) \, U(x, t) + \vartheta \, \absbig{U(x, t)}^2 \, U(x, t)\,, \\
U(x, t_0) = U_0(x)\,, \quad (x, t) \in \Omega \times [t_0, T]\,,
\end{cases}
\end{equation}  
can be cast into~\eqref{eq:GinzburgLandauType} with $\alpha_1 = \alpha \in \RR_{> 0}$, $\alpha_0 = 0$, $\beta_1 = \beta \in \RR$, and $\beta_2 = \vartheta \in \RR$.
\end{enumerate}
\section{Convergence analysis}
\label{sec:Convergence}
\MyParagraph{Guide line}
Our main objective in this section is the statement of a convergence result for high-order exponential operator splitting methods involving complex coefficients~\eqref{eq:Splitting} applied to the general model problem~\eqref{eq:GeneralModel1}.
As special cases, high-order reaction-diffusion equations~\eqref{eq:ReactionDiffusion}, complex Ginzburg--Landau-type equations~\eqref{eq:GinzburgLandauType}, and reductions to related linear equations are included in the analysis. 
For this purpose, we reconsider the approach exploited in~\cite{Thalhammer2008,Thalhammer2012} for standard real splitting methods applied to linear and nonlinear evolution equations of Schr{\"o}dinger type and indicate suitable modifications. 
In Section~\ref{sec:Stability}, we examine potential stability issues, and in Section~\ref{sec:LocalError}, we outline the derivation of local error expansions such that iterated commutators characterise the regularity requirements.
For the statement of Theorem~\ref{thm:Theorem} within the context of abstract evolution equations~\eqref{eq:IVP}, we employ the frameworks of analytic semigroups or strongly continuous unitary groups, respectively.
In view of practical implementations, the treatment of full discretisations based on time-splitting methods combined with spectral space discretisations relies on more concrete settings of separable Hilbert spaces.
Numerical experiments that confirm the validity of stability conditions and global errors bounds are presented in Section~\ref{sec:NumericalExperiments}.
They in particular imply that complex exponential operator splitting methods retain their classical orders for high-order reaction-diffusion equations and complex Ginzburg--Landau-type equations.

\MyParagraph{Fundamental means}
We here refrain from a survey of the theoretical foundations and instead refer to the standard monographs~\cite{EngelNagel2000,HillePhillips1974,Lunardi2013,Pazy1983}.
Brief descriptions of the powerful theory of sectorial operators generating analytic semigroups on Banach spaces and applications to nonautonomous as well as quasilinear parabolic problems are also found in our former works~\cite{BlanesCasasThalhammer2018,GonzalezThalhammer2006}.
\subsection{Stability conditions}
\label{sec:Stability}
\MyParagraph{Derivation}
With regard to a stability analysis of complex exponential operator splitting methods~\eqref{eq:Splitting} applied to the general model problem~\eqref{eq:GeneralModel1}, we study the decisive contribution involving the highest power of the Laplace operator and a complex constant
\begin{equation*}
\begin{cases}
v'(t) = A_K \, v(t)\,, \quad A_K = \alpha_K \, \Delta^K\,, \\
v(t_n) \text{ given}\,, \quad t \in [t_n, t_n + \taun]\,, \quad \taun > 0\,.
\end{cases} 
\end{equation*}
We recall our assumption of a spatial domain that is defined by the Cartesian product of bounded intervals and the imposed periodic boundary conditions. 
Without loss of generality, we may consider $\Omega = [- \, \pi, \pi]^d$.
A natural choice for the underlying function space is provided by the Lebesgue space of complex-valued square-integrable functions
\begin{equation*}
(L^2(\Omega, \CC), \scp{L^2}{\cdot}{\cdot}, \norm{L^2}{\cdot})\,,
\end{equation*}
endowed with the standard inner product and the induced norm.
A complete orthonormal system of this Hilbert space is given by the family of Fourier functions $(\nF_m)_{m \in \ZZ^d}$, which forms a set of eigenfunctions such that 
\begin{equation}
\label{eq:Eigenfunctions}
\begin{gathered}
\nF_m: \CC^d \longrightarrow \CC: x = (x_1, \dots, x_d) \longmapsto \tfrac{1}{(2 \, \pi)^d} \prod_{\ell=1}^{d} \ee^{\, \ii \, m_{\ell} \, (x_{\ell} + \pi)}\,, \\    
\Delta^K \, \nF_m = \big(- \, \lambda_m\big)^K \nF_m\,, \quad 
\lambda_m = \abs{m}^2 = \sum_{\ell=1}^{d} m_{\ell}^2 \geq 0\,, \\
m = (m_1, \dots, m_d) \in \ZZ^d\,. 
\end{gathered}
\end{equation}
Due to the fact that $\lambda_m \to \infty$ whenever a component of the multi-index $m \in \ZZ^d$ tends to infinity, the term $(- \, \lambda_m)^K$ associated with the highest power of the Laplace operator indeed dominates lower-order contributions.
By Fourier series representations, we obtain the relations 
\begin{equation*}
\begin{gathered}
v(t_n) = \sum_{m \in \ZZ^d} v_m(t_n) \, \nF_m\,, \\
\nE_{\taun, a_j A_K} \, v(t_n) = \sum_{m \in \ZZ^d} \ee^{\taun \, a_j \, \alpha_K \, (- \, \lambda_m)^K} v_m(t_n) \, \nF_m\,, \quad
j \in \{1, \dots, s\}\,.
\end{gathered}
\end{equation*}
Moreover, Parseval's identity implies 
\begin{equation*}
\begin{split}
&\normbig{L^2}{\nE_{\taun, a_j A_K} \, v(t_n)}^2 = \sum_{m \in \ZZ^d} \absBig{\ee^{\taun \, a_j \, \alpha_K \, (- \, \lambda_m)^K}}^2 \, \absbig{v_m(t_n)}^2\,, \\
&\qquad = \sum_{m \in \ZZ^d} \ee^{\, 2 \, (- \, 1)^K \, \Re(a_j \, \alpha_K) \, \taun \, \lambda_m^K} \, \absbig{v_m(t_n)}^2\,, \quad 
j \in \{1, \dots, s\}\,.
\end{split}
\end{equation*}
As a consequence, for positive time stepsizes $\taun > 0$, boundedness is ensured provided that
\begin{equation*}
\ee^{\, (- \, 1)^K \, \Re(a_j \, \alpha_K) \, \taun \, \lambda_m^K} \leq 1\,, \quad m \in \ZZ^d\,.
\end{equation*}
Altogether, for evolution equations involving complex constants and exponential operator splitting methods with complex coefficients, we obtain the stability conditions 
\begin{equation}
\label{eq:StabilityConditions}  
\begin{gathered}
(- \, 1)^K \, \Re\big(a_j \, \alpha_K\big) = (- \, 1)^K \, \Big(\Re(a_j) \, \Re(\alpha_K) - \Im(a_j) \, \Im(\alpha_K)\Big) \leq 0\,, \\
j \in \{1, \dots, s\}\,.
\end{gathered}
\end{equation}
In particular, they must hold for high-order reaction-diffusion equations~\eqref{eq:ReactionDiffusion} with $K = 4$ and for complex Ginzburg--Landau-type equations~\eqref{eq:GinzburgLandauType} with $K = 1$, respectively.

\MyParagraph{Special case}
It is expedient to reconsider the stability conditions~\eqref{eq:StabilityConditions} for exponential operator splitting methods involving real coefficients $(a_j)_{j=1}^{s}$.
In this special case, the following simplifications are valid
\begin{equation*}
\begin{gathered}
\Im(a_j) = 0 \quad \Longrightarrow \quad (- \, 1)^K \, \Re(a_j) \, \Re(\alpha_K) \leq 0\,, \\
j \in \{1, \dots, s\}\,. 
\end{gathered}
\end{equation*}
On the one hand, for well-posed reaction-diffusion-type equations satisfying the relation $(- \, 1)^K \, \Re(\alpha_K) < 0$, this yields 
\begin{equation*}
a_j = \Re(a_j) \geq 0\,, \quad j \in \{1, \dots, s\}\,. 
\end{equation*}
We point out that the validity of these stability conditions can be granted for higher-order exponential operator splitting methods involving complex coefficients, whereas a second-order barrier is effective for any standard real splitting method.
On the other hand, due to $\Re(\alpha_K) = 0$, stability can be ensured for Schr{\"o}dinger equations.

\MyParagraph{Conclusion}
It is noteworthy that the sizes of the real and imaginary parts of the decisive constant $\alpha_K \in \CC$ in general affect the stability behaviour of complex exponential operator splitting methods.
This conclusion from~\eqref{eq:StabilityConditions} is illustrated by a numerical experiment in Section~\ref{sec:NumericalExperiments} below, see Figure~\ref{fig:Figure6}. 
In order to avoid such a phenomenon, it is desirable to apply high-order complex splitting methods with the additional feature $(a_j)_{j=1}^{s} \in \RR^s$. 
In Table~\ref{tab:Table1}, we summarise the stability properties of various real and complex exponential operator splitting methods.
For the stated reasons, we highlight generally applicable fourth- and sixth-order schemes with non-negative coefficients $(a_j)_{j=1}^{s}$ (Methods~5, 7, 9, 10, 11, 12). 

\begin{table}[t!]
\caption{Survey of lower- and higher-order exponential operator splitting methods. 
Specification of type (real, complex), additional features (symmetric~(s), symmetric-conjugate~(sc), and alternating-conjugate~(ac)) and characteristics (nonstiff order~$p$, number of stages~$s$).
Stability properties for evolution equations of parabolic type.
For schemes involving complex coefficients, the conditions $\Re \, a_j \geq 0$ for $j \in \{1, \dots, s\}$ are necessary, but not sufficient. 
Schemes involving real coefficients $(a_j)_{j=1}^{s}$ remain stable for evolution equations of Schr{\"o}dinger type. \\[-2mm]}
{\begin{tabular}{|l|l|}\hline
Method 1 (Lie--Trotter, real, $p = s = 1$) & Stability ($a_1 > 0$) \\\hline
Method 2 (Strang, real s, $p = s = 2$) & Stability ($a_1, a_2 \geq 0$) \\\hline 
Method 3 (Yoshida, real s, $p = s = 4$) & Instability ($a_3 < 0$) \\\hline
Method 4 (Yoshida, complex s, $p = s = 4$) & $\Re \, a_1, \dots, \Re \, a_s \geq 0$ \\\hline
Method 5 (complex s, $p = 4$, $s = 6$) & Stability ($a_1, \dots, a_s \geq 0$) \\\hline
Method 6 (complex sc, $p = s = 3$) & $\Re \, a_1, \dots, \Re \, a_s \geq 0$ \\\hline
Method 7 (complex sc, $p = 3$, $s = 4$) & Stability ($a_1, \dots, a_s \geq 0$)\\\hline
Method 8 (complex sc, $p = s = 4$) & $\Re \, a_1, \dots, \Re \, a_s \geq 0$ \\\hline
Method 9 (complex sc, $p = 4$, $s = 6$) & Stability ($a_1, \dots, a_s \geq 0$) \\\hline
Method 10 (complex sc, $p = 6$, $s = 16$) & Stability ($a_1, \dots, a_s \geq 0$) \\\hline
Method 11 (complex ac, $p = 4$, $s = 7$) & Stability ($a_1, \dots, a_s \geq 0$) \\\hline
Method 12 (complex ac, $p = 6$, $s = 19$) & Stability ($a_1, \dots, a_s \geq 0$) \\\hline
\end{tabular}}
\label{tab:Table1}
\end{table}
\subsection{Local error expansions}
\label{sec:LocalError}
\MyParagraph{Fundamental means}
For the purpose of illustration, we focus on the derivation of suitable local error expansions for complex operators splitting methods~\eqref{eq:Splitting} applied to abstract evolution equations~\eqref{eq:IVP} that are defined by two linear operators.
In view of our model problem, we associate $F_1 = A$ with a differential operator and $F_2 = B$ with a multiplication operator.
Fundamental means for stepwise expansions of the exact and numerical evolution operators, which ensure the specification of the arising remainders, are provided by the variation-of-constants formula and Taylor series expansions.
The characterisation of the resulting regularity requirements is linked to the identification of iterated commutators.
The extension to the significantly more involved nonlinear case in the lines of~\cite{Thalhammer2012} relies on the formal calculus of Lie-derivatives.

\MyParagraph{Linear evolution equations}
We study linear evolution equations
\begin{equation*}
u'(t) = F\big(u(t)\big) = (A + B) \, u(t)\,, \quad t \in [t_0, T]\,,
\end{equation*}
and employ a compact notation for the associated evolution operators 
\begin{equation*}
\nE_{\tau, F} = \ee^{\, \tau (A + B)}\,, \quad \nE_{\tau, A} = \ee^{\, \tau A}\,, \quad \nE_{\tau, B} = \ee^{\, \tau B}\,.
\end{equation*}
In order to illustrate the general procedure for high-order exponential operator splitting methods, we consider a scheme of nonstiff order~$p$ involving four stages.
We meanwhile use the general form 
\begin{equation*}
\begin{gathered}
s = 4: \quad \nS_{\tau, F} = \ee^{\, B_s} \, \ee^{\, A_s} \cdots \, \ee^{\, B_1} \, \ee^{\, A_1}\,, \\
A_j = \tau \, a_j \, A\,, \quad B_j = \tau \, b_j \, B\,, \quad j \in \{1, \dots, s\}\,.
\end{gathered}
\end{equation*}
\begin{enumerate}[(i)]
\item
The derivation of appropriate expansions of the exact evolution operators is based on a repeated application of the linear variation-of-constants formula 
\begin{equation*}
u(t) = \ee^{\, (t - t_0) A} \, u(t_0) + \int_{t_0}^{t} \ee^{\, (t - \zeta) A} \, B \, u(\zeta) \; \dd \zeta\,, \quad t \in [t_0, T]\,.
\end{equation*}
Under suitable regularity restrictions, the resulting representations
\begin{equation}
\label{eq:ExpansionExact}
\nE_{\tau, F} = \ee^{\, \tau (A + B)} = \sum_{k=0}^{p} \nE^{(k)}_{\tau, F} + \nO\big(\tau^{p+1}\big)\,,  
\end{equation}
where the dominant terms are given by iterated integrals such as 
\begin{equation*}
\begin{gathered}
\nE^{(0)}_{\tau, F} = \ee^{\, \tau A}\,, \\
\nE^{(1)}_{\tau, F} = \int_{0}^{\tau} \ee^{\, (\tau - \zeta_1) A} \, B \, \ee^{\, \zeta_1 A} \; \dd \zeta_1\,, \\
\nE^{(2)}_{\tau, F} = \int_{0}^{\tau} \int_{0}^{\zeta_1} \ee^{\, (\tau - \zeta_1) A} \, B \, \ee^{\, (\zeta_1 - \zeta_2) A} \, B \, \ee^{\, \zeta_2 A} \; \dd \zeta_2 \, \dd \zeta_1\,,
\end{gathered}
\end{equation*}
are well-defined.
Detailed calculations and explanations, valid for the general case, are found in~\cite{Thalhammer2008}. 
\item
For the numerical evolution operators, the derivation of appropriate representations that resemble~\eqref{eq:ExpansionExact} are obtained by Taylor series expansions of the evolution operators associated with the multiplication operators.
This is done step-by-step such that the remainders reflect the regularity requirements in accordance with numerical evidence. 
For detailed calculations and explanations, we again refer to~\cite{Thalhammer2008}. 
In the case of the above stated example, this yields a relation of the form 
\begin{equation}
\label{eq:ExpansionNumerical}
\nS_{\tau, F} = \ee^{\, B_4} \, \ee^{\, A_4} \, \ee^{\, B_3} \, \ee^{\, A_3} \, \ee^{\, B_2} \, \ee^{\, A_2} \, \ee^{\, B_1} \, \ee^{\, A_1}
= \sum_{k=0}^{p} \nS^{(k)}_{\tau, F} + \nO\big(\tau^{p+1}\big)\,,
\end{equation}
where the leading contributions are given by
\begin{equation*}
\begin{split}
\nS^{(0)}_{\tau, F} &= \ee^{\, A_4 + A_3 + A_2 + A_1}\,, \\
\nS^{(1)}_{\tau, F}
&= B_4 \, \ee^{\, A_4 + A_3 + A_2 + A_1} + \ee^{\, A_4} \, B_3 \, \ee^{\, A_3 + A_2 + A_1} \\
&\qquad + \ee^{\, A_4 + A_3} \, B_2 \, \ee^{\, A_2 + A_1} + \ee^{\, A_4 + A_3 + A_2} \, B_1 \, \ee^{\, A_1}\,, \\
\nS^{(2)}_{\tau, F}
&= \tfrac{1}{2} \, B_4^2 \, \ee^{\, A_4 + A_3 + A_2 + A_1} + B_4 \, \ee^{\, A_4} \, B_3 \, \ee^{\, A_3 + A_2 + A_1} \\
&\qquad + B_4 \, \ee^{\, A_4 + A_3} \, B_2 \, \ee^{\, A_2 + A_1} + B_4 \, \ee^{\, A_4 + A_3 + A_2} \, B_1 \, \ee^{\, A_1} \\
&\qquad + \tfrac{1}{2} \, \ee^{\, A_4} \, B_3^2 \, \ee^{\, A_3 + A_2 + A_1} + \ee^{\, A_4} \, B_3 \, \ee^{\, A_3} \, B_2 \, \ee^{\, A_2 + A_1} \\
&\qquad + \ee^{\, A_4} \, B_3 \, \ee^{\, A_3 + A_2} \, B_1 \, \ee^{\, A_1} + \tfrac{1}{2} \, \ee^{\, A_4 + A_3} \, B_2^2 \, \ee^{\, A_2 + A_1} \\
&\qquad + \ee^{\, A_4 + A_3} \, B_2 \, \ee^{\, A_2} \, B_1 \, \ee^{\, A_1} + \tfrac{1}{2} \, \ee^{\, A_4 + A_3 + A_2} \, B_1^2 \, \ee^{\, A_1}\,.
\end{split}
\end{equation*}
\item
The order conditions for complex exponential operator splitting methods applied to linear evolution equations involving unbounded operators follow from a further analysis of the differences 
\begin{equation*}
\sum_{k=0}^{p} \Big(\nS^{(k)}_{\tau, F} - \nE^{(k)}_{\tau, F}\Big) = \nO\big(\tau^{p+1}\big)\,, 
\end{equation*}
see~\eqref{eq:ExpansionExact} and~\eqref{eq:ExpansionNumerical}.
For this pupose, the contributions in the expansions of the numerical evolution operators are understood as quadrature approximations to the integrals arising in the expansions of the exact evolution operators.
Commutators naturally result from suitable Taylor series expansions, which comprise elements such as 
\begin{equation*}
\begin{gathered}
\tfrac{\dd}{\dd \zeta} \, \ee^{\, (\tau - \zeta) A} \, B \, \ee^{\, \zeta A} = - \, \ee^{\, (\tau - \zeta) A} \, \text{ad}_A(B) \, \ee^{\, \zeta A}\,, \\
\tfrac{\dd^2}{\dd \zeta^2} \, \ee^{\, (\tau - \zeta) \, A} \, B \, \ee^{\, \zeta A} = \ee^{\, (\tau - \zeta) A} \, \text{ad}_A^{\, 2}(B) \, \ee^{\, \zeta A}\,, 
\end{gathered}
\end{equation*}
see also~\eqref{eq:IteratedCommutators}.
\end{enumerate}

\MyParagraph{Generalisation to nonlinear evolution equations}
The derivation of suitable local error expansions for complex exponential operator splitting methods applied to a nonlinear evolution equation of the form~\eqref{eq:IVP}, retained from the general model~\eqref{eq:GeneralModel1}--\eqref{eq:GeneralModel2}, follows the approach sketched above.
Important means are hereby provided by the formal calculus of Lie-derivatives.
In particular, this concerns the application of the nonlinear variation-of-constants formula to obtain appropriate expansions of the exact evolution operators and the occurrence of iterated commutators involving nonlinear operators, see~\eqref{eq:IteratedCommutators}.

\MyParagraph{Regularity requirements}
As detailed in~\cite{Thalhammer2008}, a fundamental ingredient in the rigorous treatment of high-order exponential operator splitting methods is the specification of the remainders and the characterisation of the arising commutators. 
In addition, the resulting regularity requirements are related to fractional power spaces of sectorial operators or, more concretely, to Sobolev spaces.
For the general model problem~\eqref{eq:GeneralModel1}--\eqref{eq:GeneralModel2}, where $F_1$ represents a linear differential operator and~$F_2$ a nonlinear multiplication operator, this implies that the nonstiff orders $p \in \NN_{\geq 1}$ of the considered complex exponential operator splitting methods~\eqref{eq:Splitting} are retained, provided that the problem data and hence the solution is sufficiently regular. 
More precisely, we can draw the following conclusions. 
\begin{enumerate}[(i)]
\item   
For a linear evolution equation involving the Laplace operator and a regular space-dependent function~$W$, straightforward calculations show that the commutator $\text{ad}_{\Delta}(W)$ defines a differential operator of order one that comprises derivatives of~$W$ up to order two.
By induction, it is seen that $\text{ad}_{\Delta}^{\, p}(W)$ leads to a differential operator of order~$p$ that involves derivatives of~$W$ up to order~$2 \, p$.
As a consequence, we retain the nonstiff order of a complex exponential operator splitting method, if the prescribed initial state and thus the exact solution is contained in the domain
\begin{equation*}
D = D\big((- \, \Delta)^{p/2}\big)\,.
\end{equation*}
\item   
In contrast, for nonlinear problems such as complex Ginzburg--Landau-type equations, we cannot expect that certain derivatives cancel, that is, the commutator $\text{ad}_{F_1}(F_2)$ corresponds to a differential operator of order~$2 \, p$. 
Hence, the nonstiff orders of convergence are retained, whenever the exact solution values remain bounded in
\begin{equation*}
D = D\big((- \, \Delta)^p\big)\,.
\end{equation*}
\end{enumerate}
Similar arguments are valid for high-order reaction-diffusion equations and associated linear equations.

\MyParagraph{Order reductions}
In case the well-definedness of $\text{ad}_{F_1}^{\, \ell}(F_2)$ is ensured for $\ell < p$ only, order reductions are encountered.
An illustrative numerical experiment for higher-order real splitting methods applied to the Gross--Pitaevskii equation is given in~\cite[Figure~4.2]{Thalhammer2012}.
\subsection{Global error estimates}
\label{sec:Theorem}
The above considerations and exemplifications lead us to the following convergence result for complex exponential operator splitting methods~\eqref{eq:Splitting}. 

\MyParagraph{Setting}
The scope of applications includes nonlinear partial differential equations of parabolic and Schr{\"o}dinger type that can be cast into the general form~\eqref{eq:GeneralModel1}.
On the one hand, for high-order reaction-diffusion equations~\eqref{eq:ReactionDiffusion}, complex Ginzburg--Landau equations~\eqref{eq:CGL}, and parabolic counterparts of Gross--Pitaevskii equations~\eqref{eq:GPEParabolic}, we employ the framework of analytical semigroups generated by sectorial operators on Banach spaces.
It is noteworthy that this framework permits to include additional lower-order partial derivatives in the linear differential operators, when desired.
On the other hand, for Gross--Pitaevskii equations~\eqref{eq:GPE}, we instead employ the framework of strongly continuous unitary groups. 

\MyParagraph{Telescopic identity}
As indicated in Section~\ref{sec:SplittingMethods}, a standard telescopic identity reveals that estimates for the differences between numerical and exact solution values 
\begin{equation*}
\begin{split}
&u_n - u(t_n) = \big(\nS_{\tau_{n-1}, F} \circ \dots \circ \nS_{\tau_0, F}\big)(u_0)  - \big(\nE_{\tau_{n-1}, F} \circ \dots \circ \nE_{\tau_0, F}\big)\big(u(t_0)\big) \\ 
&= \big(\nS_{\tau_{n-1}, F} \circ \dots \circ \nS_{\tau_0, F}\big)(u_0) - \big(\nS_{\tau_{n-1}, F} \circ \dots \circ \nS_{\tau_0, F}\big)\big(u(t_0)\big) \\
&\qquad + \sum_{\ell=1}^n \bigg(\big(\nS_{\tau_{n-1}, F} \circ \dots \circ \nS_{\tau_{\ell}, F} \circ \nS_{\tau_{\ell-1}, F}\big)\big(u(t_{\ell-1})\big) \\
&\qquad\qquad\,\, - \big(\nS_{\tau_{n-1}, F} \circ \dots \circ \nS_{\tau_{\ell}, F} \circ \nE_{\tau_{\ell-1}, F}\big)\big(u(t_{\ell-1})\big)\,, \quad n \in \{0, 1, \dots, N\}\,,
\end{split}
\end{equation*}
can be reduced to stability and local error bounds such as 
\begin{equation*}
\begin{gathered}
\normbig{X}{\big(\nS_{\tau_{n-1}, F} \circ \dots \circ \nS_{\tau_0, F}\big)(v) - \big(\nS_{\tau_{n-1}, F} \circ \dots \circ \nS_{\tau_0, F}\big)(w)} \leq C \, \norm{X}{v - w}\,, \\
\normbig{X}{\nS_{\tau_{n-1}, F}(v) - \nE_{\tau_{n-1}, F}(v)} \leq C \, \tau_{n-1}^{p+1}\,, \\ 
n \in \{0, 1, \dots, N\}\,, 
\end{gathered}
\end{equation*}
see also Lemma~5 and Lemma~9 in~\cite{Thalhammer2012}. 
Based on the stability analysis of Section~\ref{sec:Stability} and the derivation of local error estimates as sketched in Section~\ref{sec:LocalError}, we thus obtain the following statement. 
We recall the reformulation of the general model problem~\eqref{eq:GeneralModel1}--\eqref{eq:GeneralModel2} as abstract evolution equation~\eqref{eq:IVP}, where~$F_1$ represents a linear differential operator that is given by linear combinations of powers of the Laplacian and~$F_2$ a nonlinear multiplication operator.

\begin{Theorem}
\label{thm:Theorem}
Let $(X, \norm{X}{\cdot})$ denote the underlying Banach space.
Consider the nonlinear evolution equation~\eqref{eq:IVP} involving the (unbounded) operators $F_{\ell}: D(F_{\ell}) \subseteq X \to X$, $\ell \in \{1, 2\}$. 
Assume that~$F_1$ generates an analytic semigroup $(\nE_{t, F_1})_{t \in \RR_{\geq 0}}$ or a strongly continuous unitary group $(\nE_{t, F_1})_{t \in \RR}$, respectively.
Suppose that the coefficients of the complex exponential operator splitting method~\eqref{eq:Splitting} fulfill the classical order conditions for some integer $p \in \NN_{\geq 1}$ and that in particular the validity of the stability bounds
\begin{equation*}
\normbig{X \leftarrow X}{\nE_{t, a_j F_1}} \leq \ee^{\, C_1 t}\,, \quad t \in [t_0, T]\,, \quad j \in \{1, \dots, s\}\,,
\end{equation*}
is ensured.
Then, provided that the exact solution values are bounded
\begin{equation*}
\normbig{D}{u(t)} \leq C_2\,, \quad t \in [t_0, T]\,, 
\end{equation*}
with respect to the norm of a suitably restricted subspace $D \subseteq X$ that is defined by the requirement that the iterated commutators arising in the expansion of the local error remain bounded
\begin{equation*}
\normbig{X \leftarrow D}{\text{ad}_{F_1}^{\, \ell}(F_2)} \leq C_3\,, \quad \ell \in \{0, 1, \dots, p\}\,,
\end{equation*}
the following global error estimate holds
\begin{equation*}
\begin{gathered}  
\normbig{X}{u_n - u(t_n)} \leq C \, \Big(\normbig{X}{u_0 - u(t_0)} + \tau_{\max}^p\Big)\,, \\
n \in \{1, \dots, N\}\,, \quad \tau_{\max} = \max_{n \in \{0, 1, \dots, N-1\}} \taun\,. 
\end{gathered}  
\end{equation*}
The positive constant $C > 0$ depends on $C_1, C_2, C_3 > 0$ and the final time, but is independent of the time increment and the number of time steps. 
\end{Theorem}

\MyParagraph{Full discretisations}
We point out that the analogue to Theorem 3~in~\cite{Thalhammer2012} for space-time discretisations by high-order time-splitting Fourier pseudo-spectral methods holds as well.
Its derivation relies on the setting of self-adjoint operators on Hilbert spaces.
Specifically, the knowledge of a complete orthonormal system of eigenfunctions associated with the Laplacian is utilised, see~\eqref{eq:Eigenfunctions}.


\begin{figure}[t!]
\begin{center}
\includegraphics[scale=0.6]{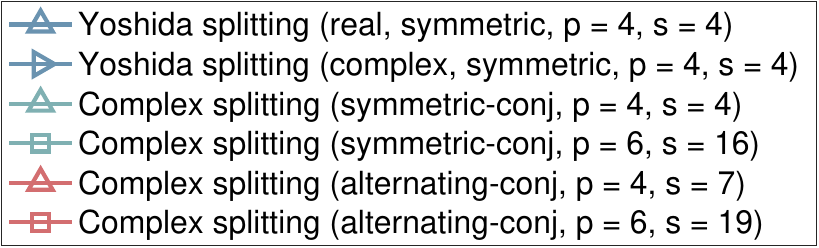}
\end{center}
\caption{Fourth- and sixth-order exponential operator splitting methods applied in numerical experiments, see Methods~3, 4, 8, 10, 11, 12 in Table~\ref{tab:Table1}.}
\label{fig:Figure1}
\end{figure}

\begin{figure}[t!]
\begin{center}
\includegraphics[scale=0.36]{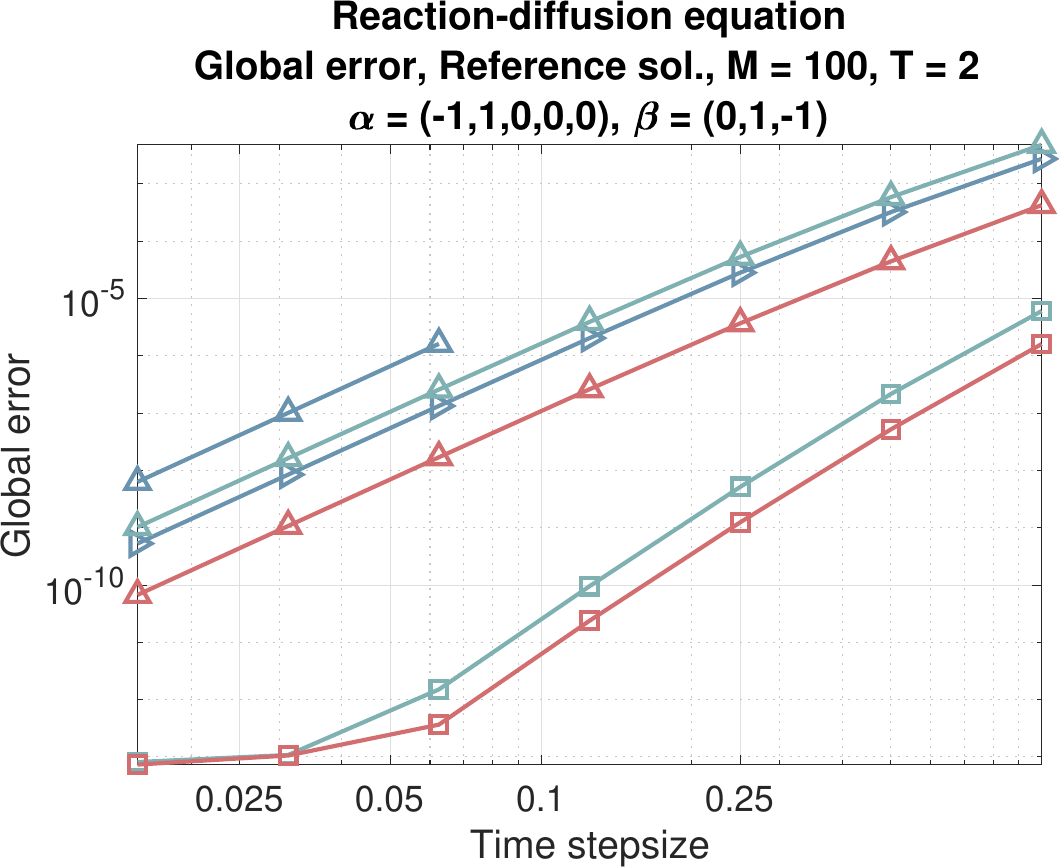} \quad 
\includegraphics[scale=0.36]{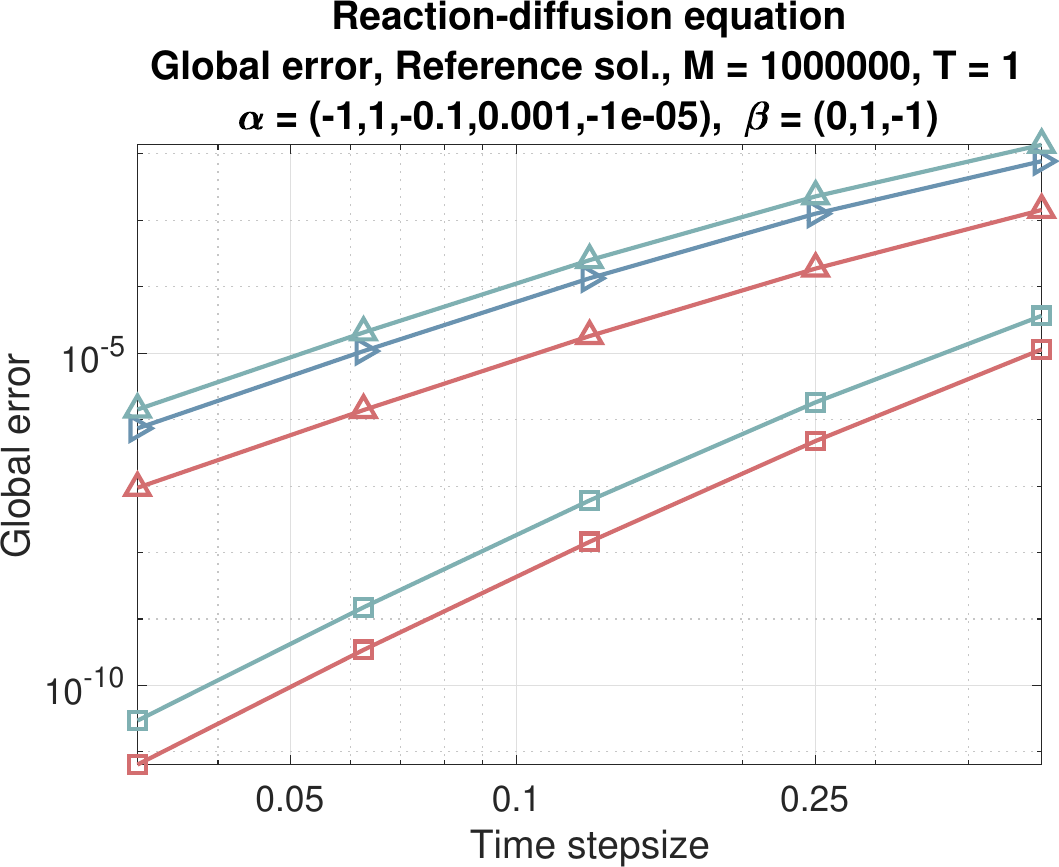} \\[2mm]
\includegraphics[scale=0.36]{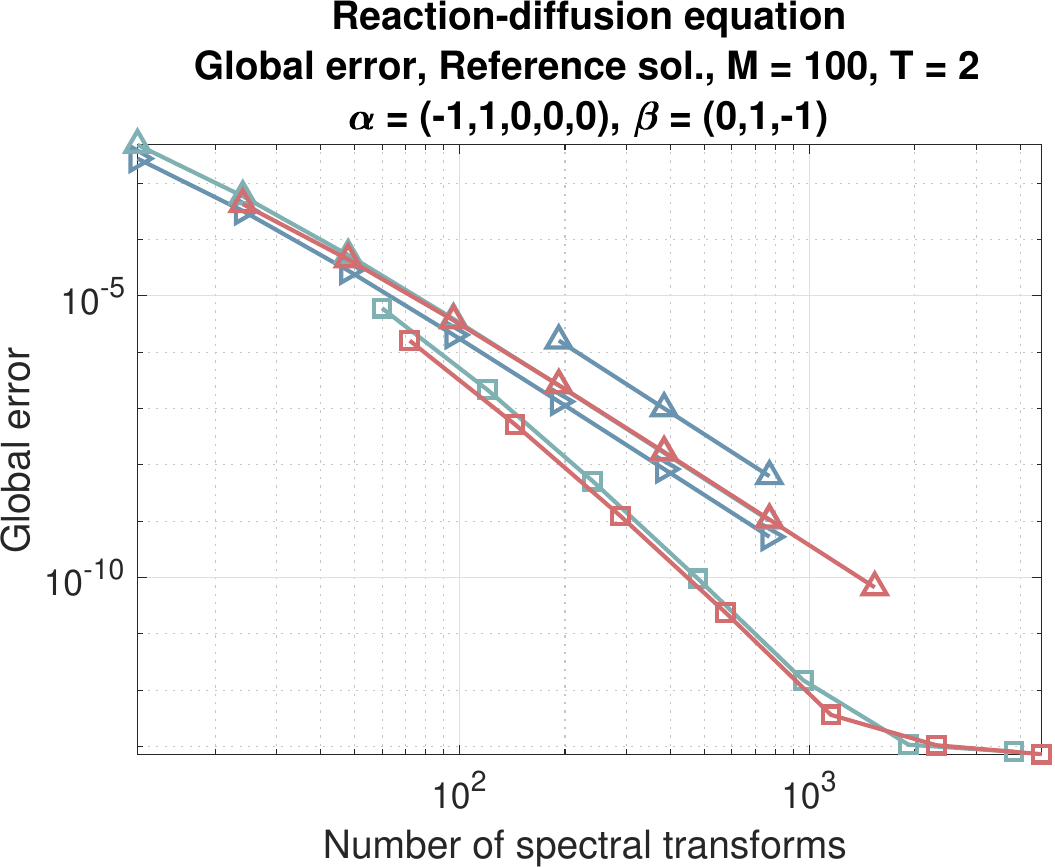} \quad 
\includegraphics[scale=0.36]{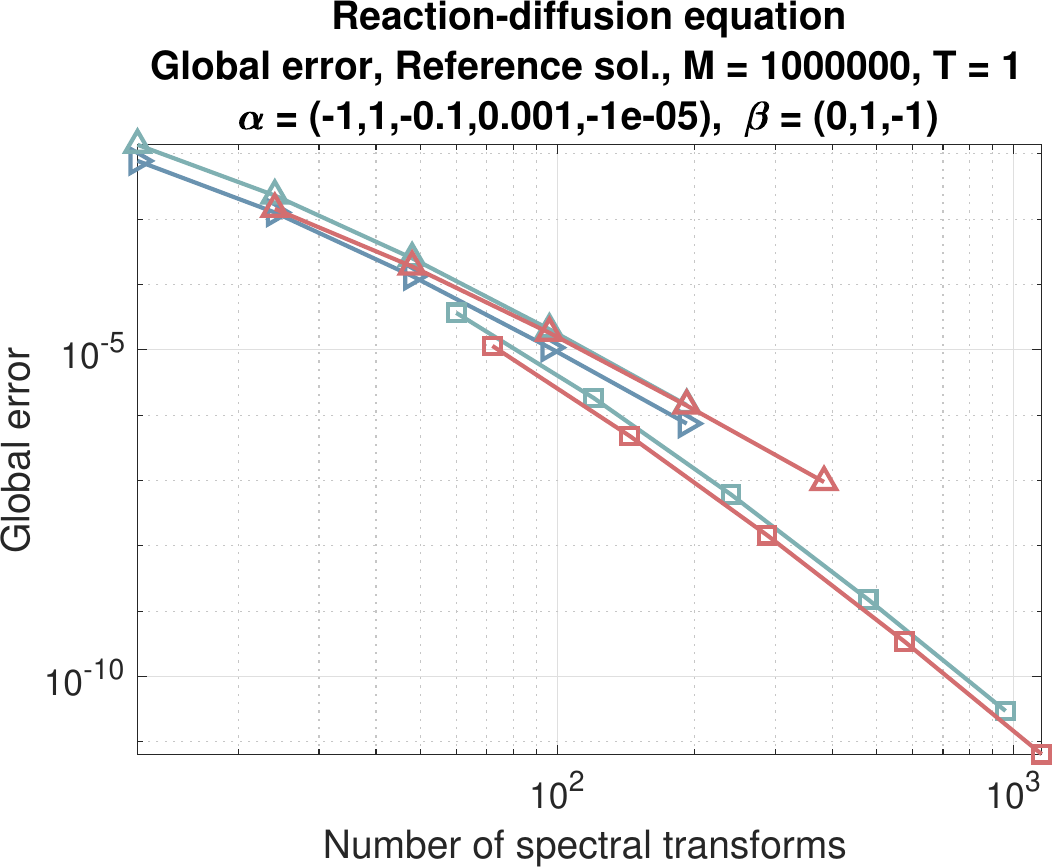} 
\end{center}
\caption{Short-term integrations of a one-dimensional reaction-diffusion equation (left) and a three-dimensional high-order reaction-diffusion equation (right) by real and complex fourth- and sixth-order exponential operator splitting methods, see Figure~\ref{fig:Figure1}.
Global errors versus time stepsizes (first row) and total numbers of spectral transforms (second row), respectively.
The complex exponential operator splitting methods remain stable and retain their classical orders, whereas the real splitting method by Yoshida suffers from severe instabilities.}
\label{fig:Figure2}
\end{figure}

\begin{figure}[t!]
\begin{center}
\includegraphics[scale=0.36]{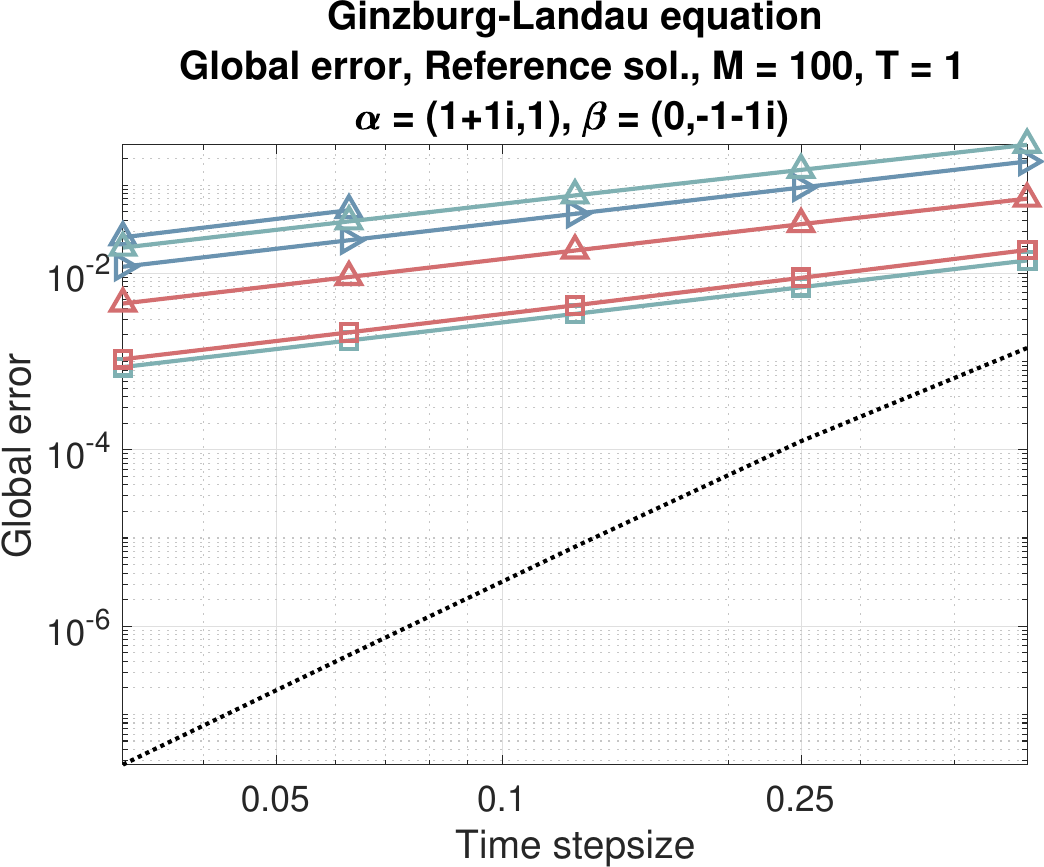} \quad 
\includegraphics[scale=0.36]{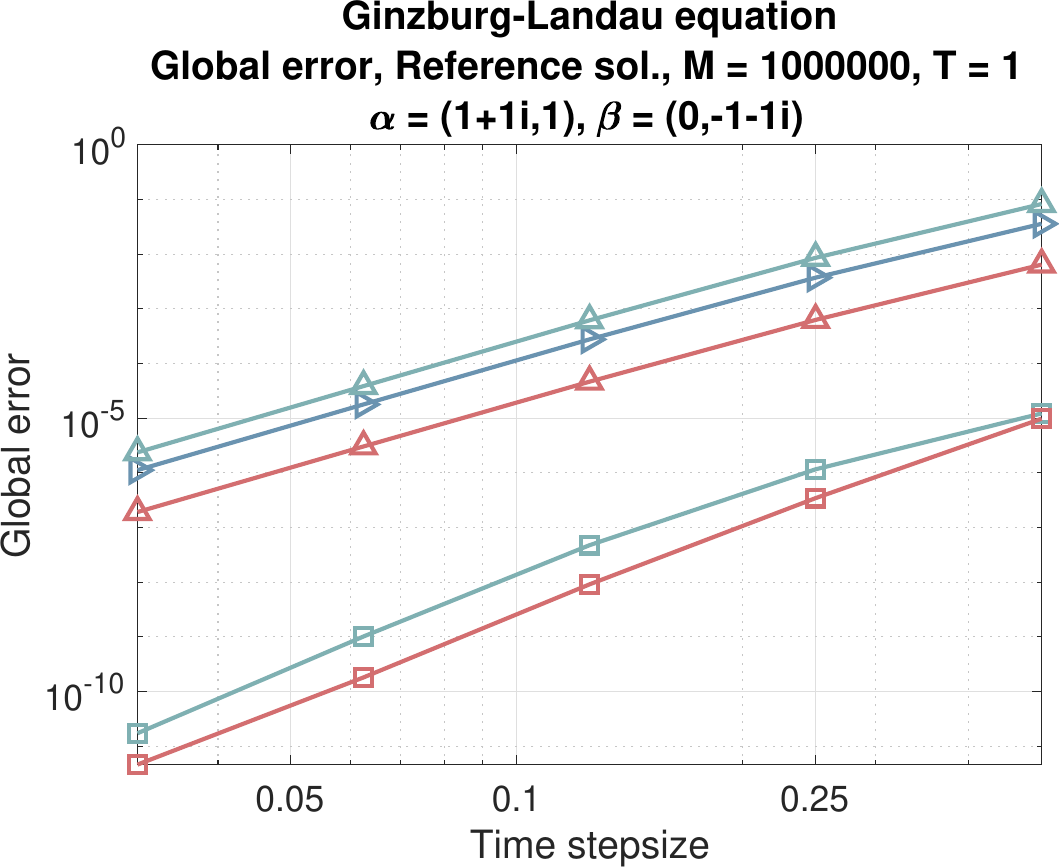} 
\end{center}
\caption{Short-term integration of complex Ginzburg--Landau equations by real and complex fourth- and sixth-order exponential operator splitting methods, see Figure~\ref{fig:Figure1}.
Global errors versus time stepsizes.
Left: Standard implementation for a one-dimensional problem. 
The complex exponential operator splitting methods remain stable, whereas the real splitting method by Yoshida becomes unstable for larger time stepsizes.
The black reference line corresponds to slope four and shows that all real and complex schemes suffer from severe order reductions.
Right: Correct implementation for a related three-dimensional problem.
All complex exponential operator splitting methods remain stable and retain their nonstiff orders. 
}
\label{fig:Figure4}
\end{figure}

\begin{figure}[t!]
\begin{center}
\includegraphics[scale=0.36]{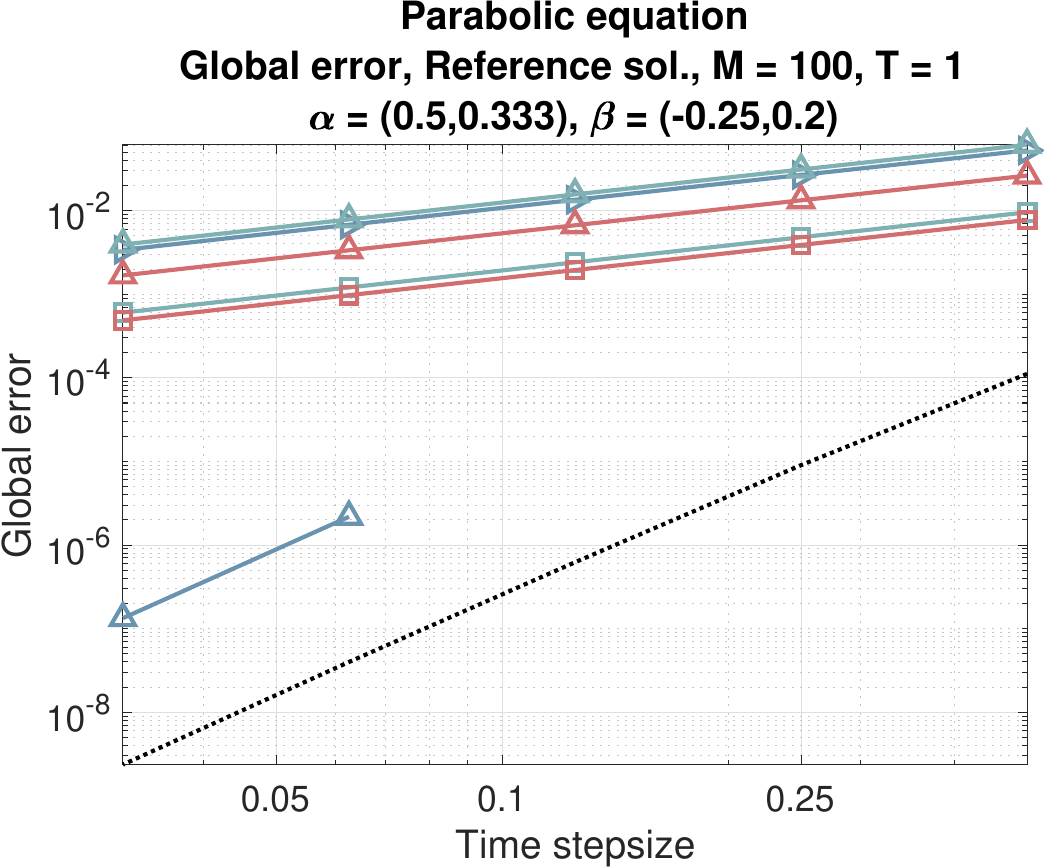} \quad 
\includegraphics[scale=0.36]{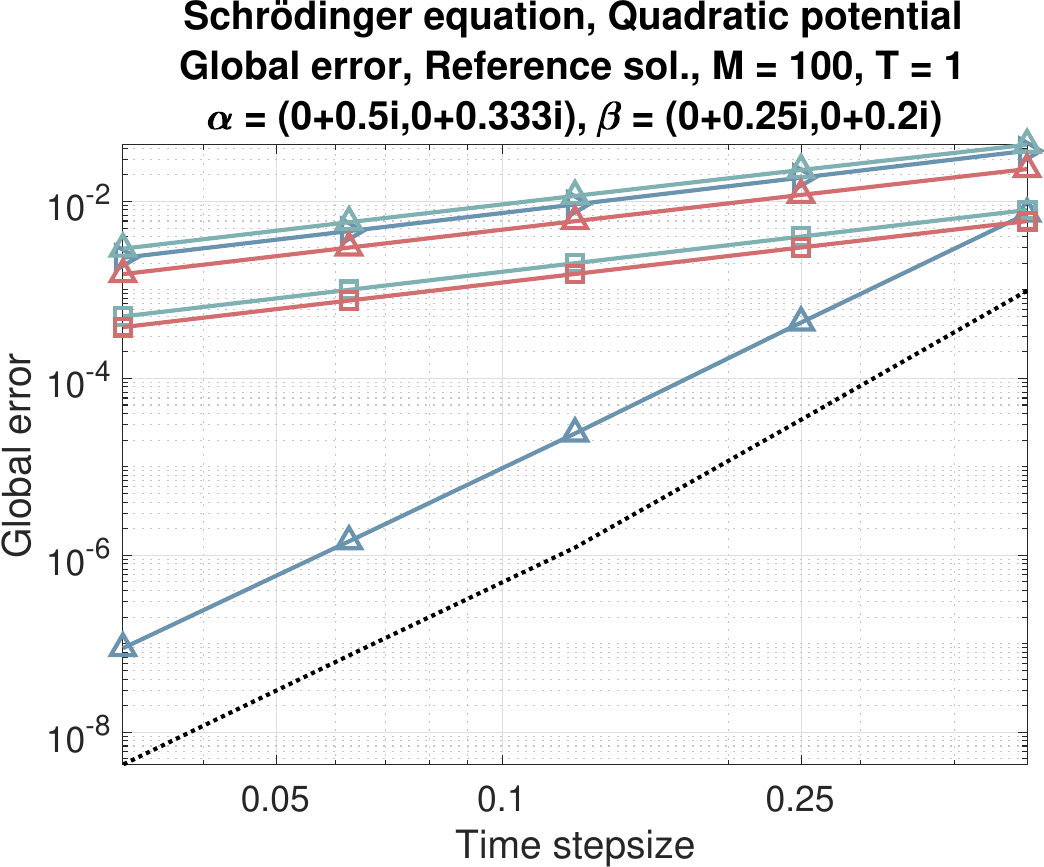} 
\end{center}
\caption{Corresponding results for related equations of parabolic and Schr{\"o}dinger type.
For a standard implementation of higher-order complex exponential operator splitting methods, severe order reductions are observed. 
}
\label{fig:Figure5}
\end{figure}

\begin{figure}[t!]
\begin{center}
\includegraphics[scale=0.36]{Fig_GlobalError_18_3d} \quad 
\includegraphics[scale=0.36]{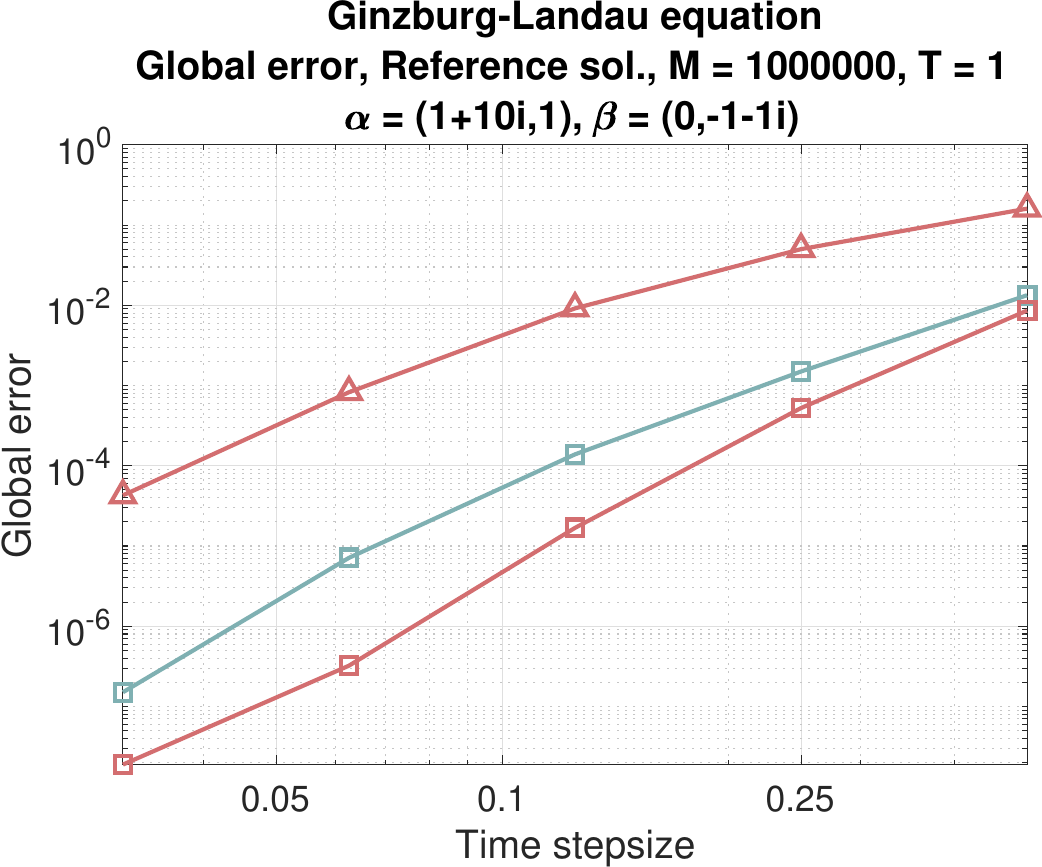} \\[2mm]
\includegraphics[scale=0.36]{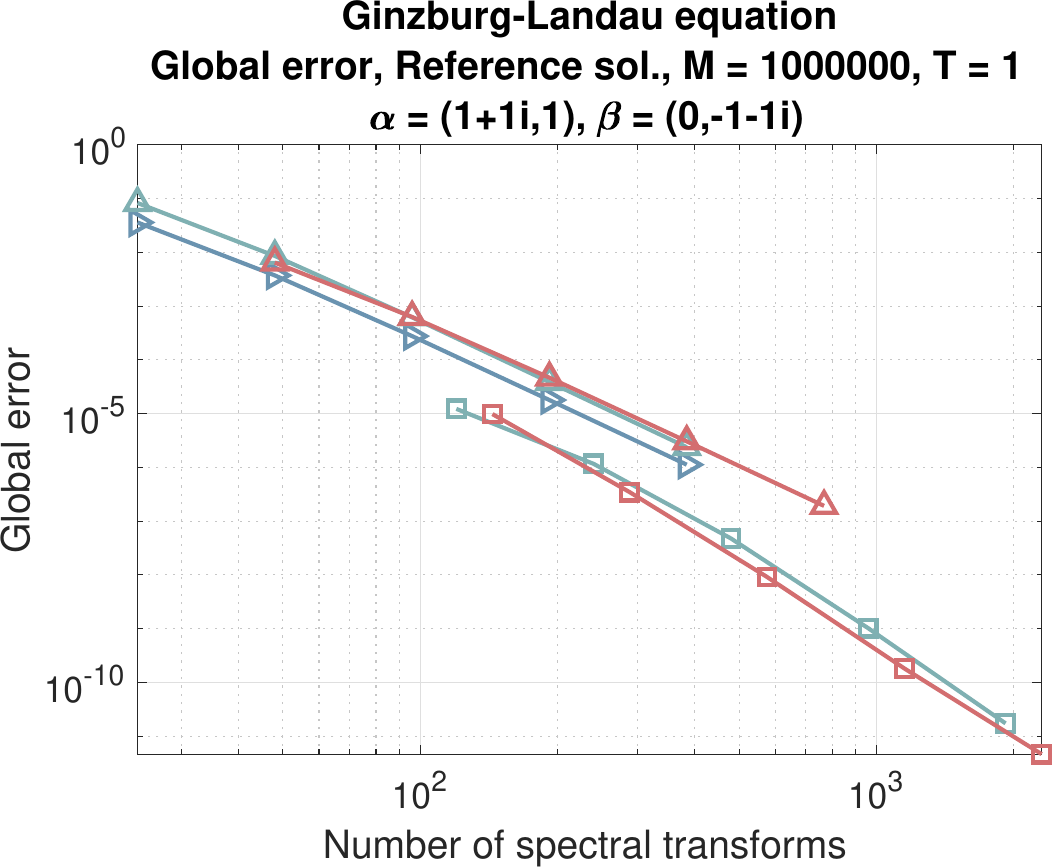} \quad 
\includegraphics[scale=0.36]{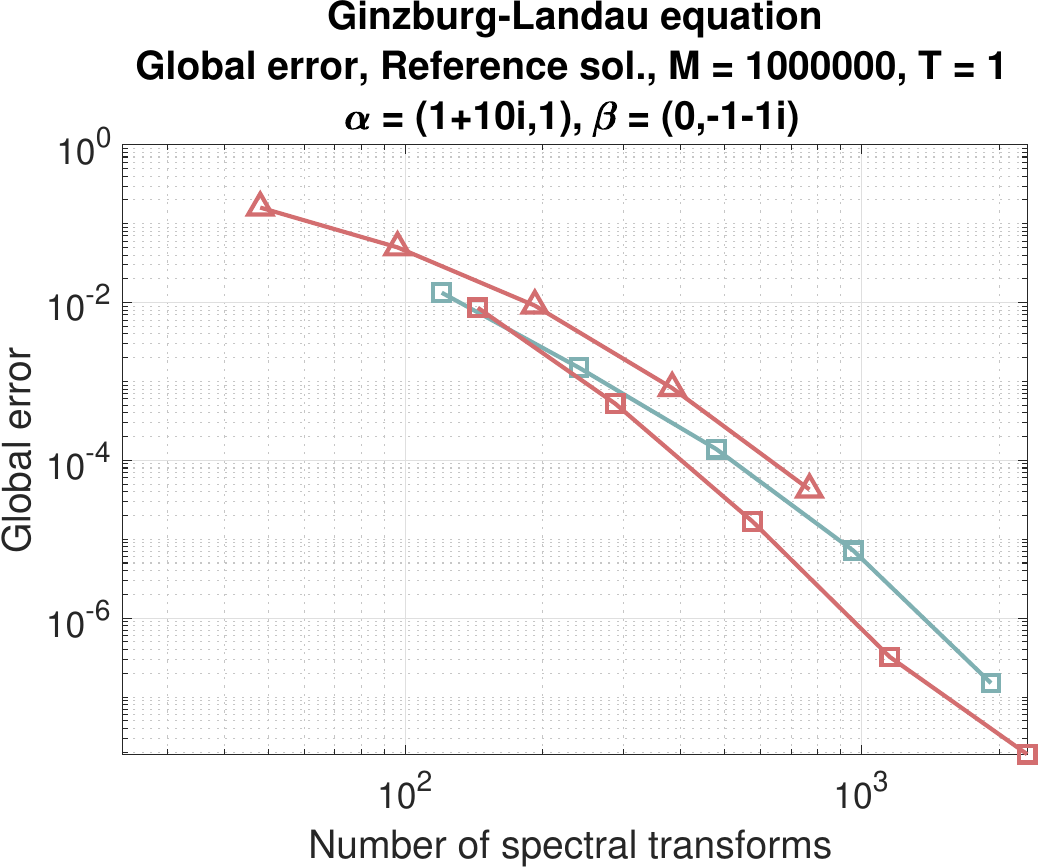} 
\end{center}
\caption{Short-term integrations of three-dimensional complex Ginzburg--Landau equations by real and complex fourth- and sixth-order exponential operator splitting methods, see Figure~\ref{fig:Figure1}.
Global errors versus time stepsizes (first row) and total numbers of spectral transforms (second row), respectively.
For $\alpha_1 = 1 + \ii$, all complex exponential operator splitting methods remain stable and retain their classical orders.
For $\alpha_1 = 1 + 10 \, \ii$, the stability conditions~\eqref{eq:StabilityConditions} have an impact on the schemes with complex coefficients $(a_j)_{j=1}^s$, whereas the schemes involving non-negative coefficients $(a_j)_{j=1}^s$ maintain stability.}
\label{fig:Figure6}
\end{figure}


\begin{figure}[t!]
\begin{center}
\includegraphics[scale=0.53]{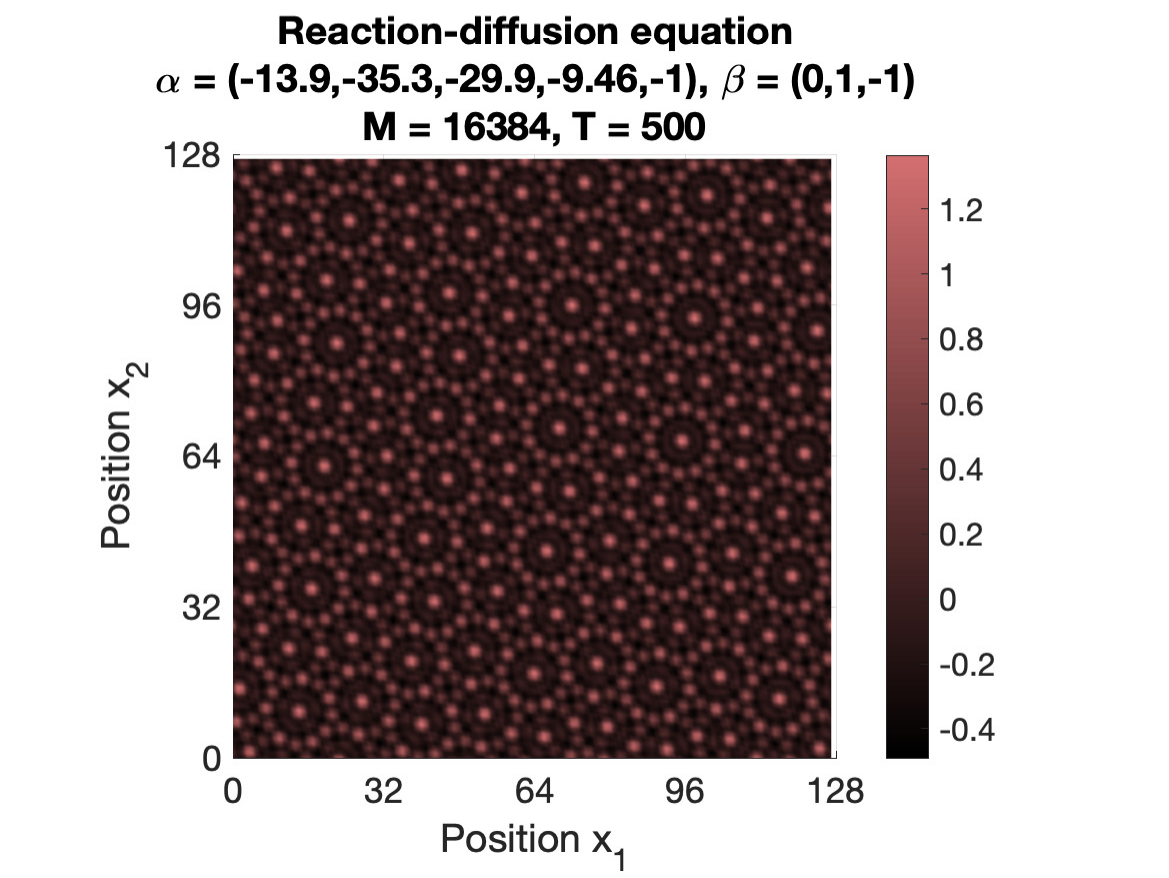} \\[2mm]
\includegraphics[scale=0.53]{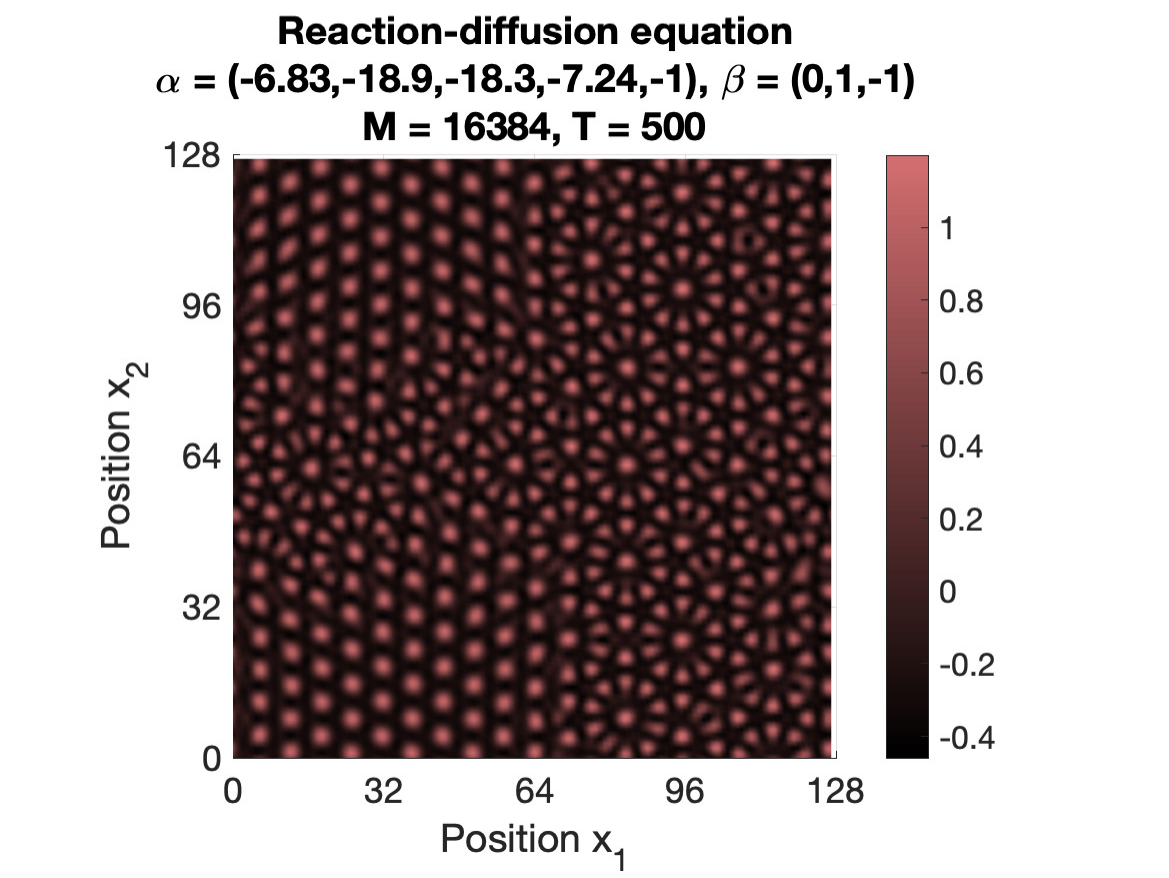} 
\end{center}
\caption{Long-term integrations of two-dimensional high-order reaction-diffusion equations forming quasicrystalline patterns.
Solution profiles at time $T = 500$.  
}
\label{fig:Figure3}
\end{figure}

\begin{figure}[t!]
\begin{center}
\includegraphics[scale=0.53]{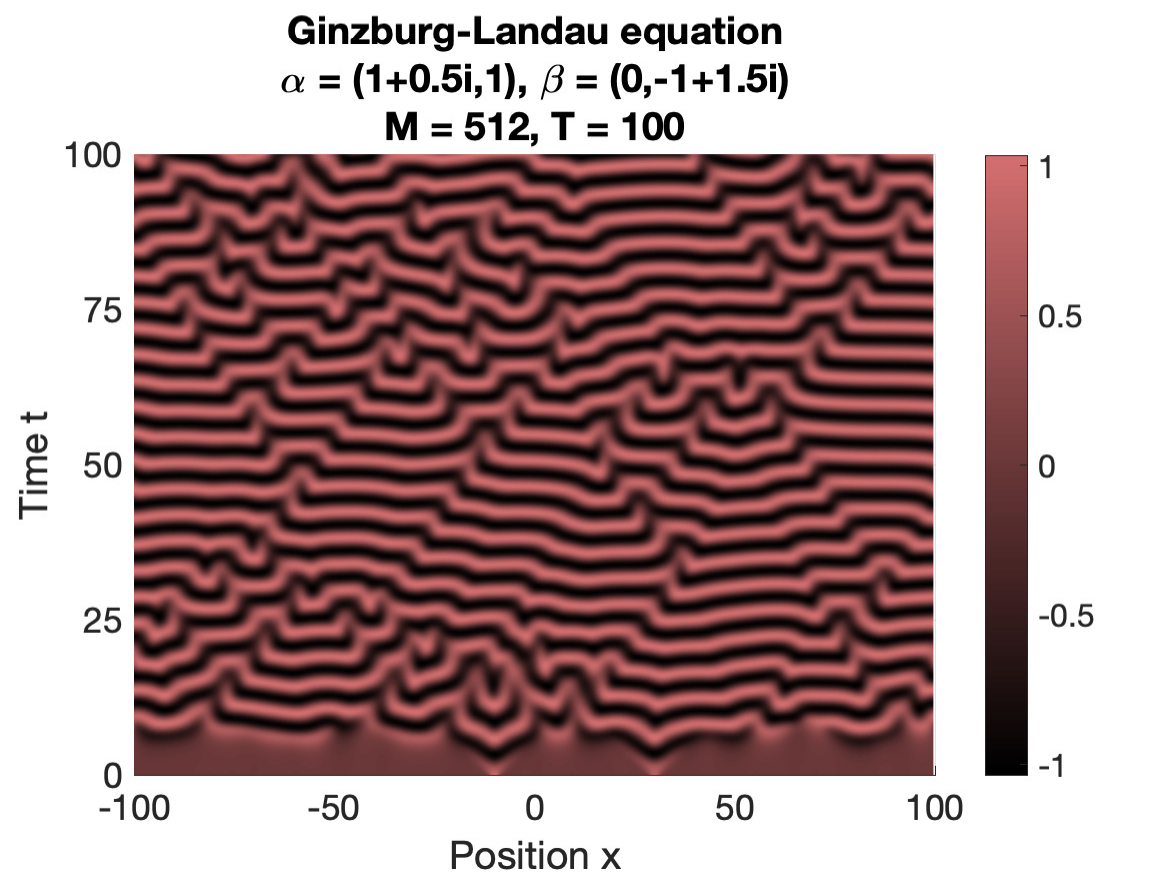} \\[2mm]
\includegraphics[scale=0.53]{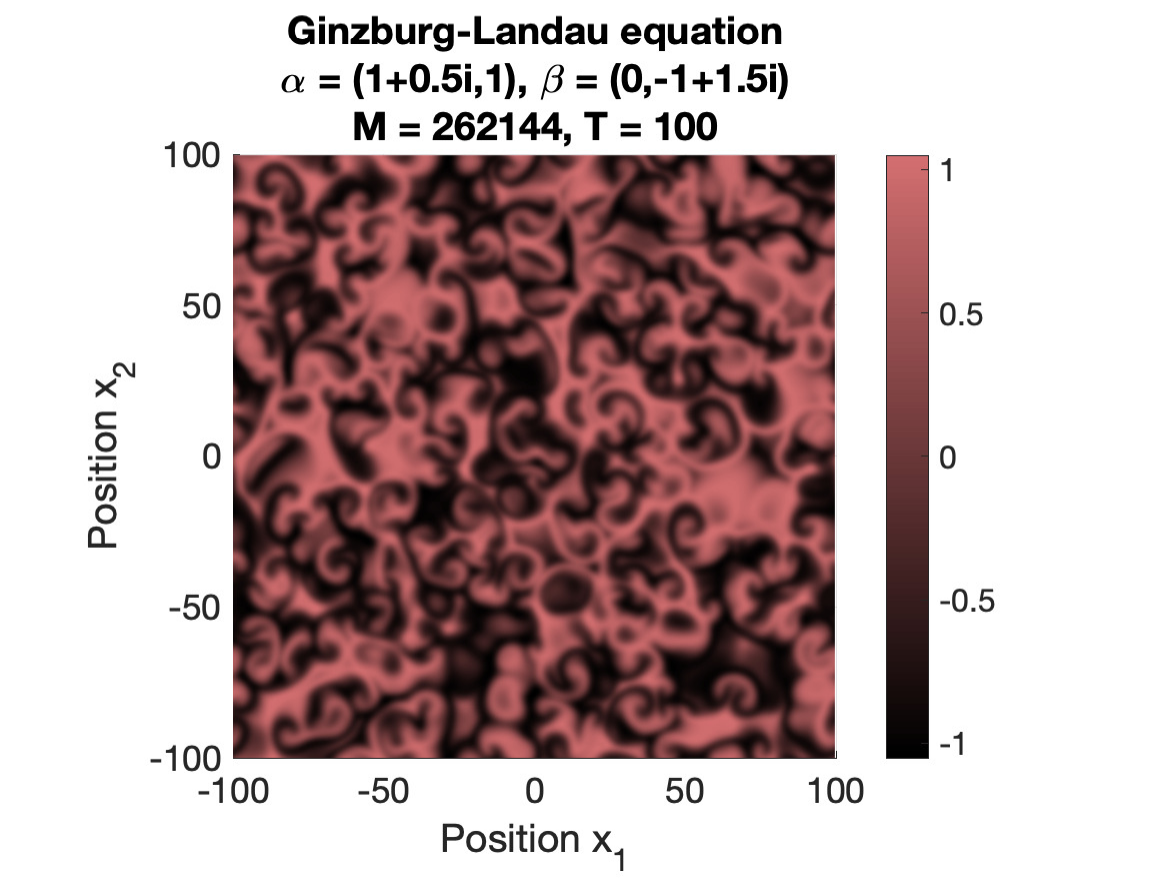} 
\end{center}
\caption{Long-term integrations of one-dimensional (up) and two-dimensional (down) complex Ginzburg--Landau equations modelling nonlinear waves.
Solution profiles at time $T = 100$.  
}
\label{fig:Figure7}
\end{figure}

\section{Numerical experiments}
\label{sec:NumericalExperiments}
In the following, we demonstrate the favourable performance of higher-order complex exponential operator splitting methods in short- and long-term integrations. 
We in particular confirm the validity of the stability conditions~\eqref{eq:StabilityConditions} and Theorem~\ref{thm:Theorem}, which implies that complex exponential operator splitting methods retain their classical orders for high-order reaction-diffusion equations and complex Ginzburg--Landau-type equations.

\MyParagraph{Source}
A \textsc{Matlab} code that illustrates the implementation of exponential operator splitting methods for high-order reaction-diffusion equations and in particular reproduces the results shown in Figure~\ref{fig:Figure2} is available at 
\begin{center}
\url{doi.org/10.5281/zenodo.13834638}.
\end{center}
The coefficients of the real and complex schemes listed in Table~\ref{tab:Table1} are found at the end of the file. 

\MyParagraph{Splitting methods}
As characteristic instances, we select the real and complex fourth-order symmetric splitting methods proposed by Yoshida~\cite{Yoshida1990} as well as two symmetric-conjugate and two alternating-conjugate splitting methods, in each case a fourth-order and a sixth-order scheme, see Figure~\ref{fig:Figure1}.

\MyParagraph{Reaction-diffusion equations}
In Figure~\ref{fig:Figure2}, we include the numerical results obtained for two test cases, a one-dimensional reaction-diffusion equation as well as a computationally more demanding high-order reaction-diffusion equation in three space dimensions.
As initial state, we choose a localised and highly regular Gaussian-like function.  
The quantity~$M$ captures the total numbers of Fourier basis functions.  
We compute a numerical reference solution and the errors at the final time with respect to the Euclidean norm.  
On the one hand, we display the global errors versus the time stepsizes.
Thus, the slopes of the lines reflect the orders of the splitting methods. 
On the other hand, we display the global errors versus the total numbers of spectral transforms (FFT, IFFT), which corresponds to the main costs of the computations and hence relates to the efficiency of the splitting methods.  
In both situations, the complex exponential operator splitting methods remain stable and retain their classical orders.
In contrast, due to the presence of negative coefficients, the real splitting method by Yoshida suffers from severe instabilities and yields no output for the three-dimensional high-order reaction-diffusion equation.
Altogether, we observe that the sixth-order schemes are superior in efficiency, even when lower tolerances are desirable. 

\MyParagraph{Ginzburg--Landau-type equations}
We point out that the realisation of higher-order exponential operator splitting methods for complex Ginzburg--Landau-type equations is a subtle issue.
Due to the presence of non-analytic nonlinearities and complex constants, standard implementations would lead to significant order reductions.
Instead, an appropriate reformulation of~\eqref{eq:GinzburgLandauType} as a system for the solution and its complex conjugate is required. 
Additional considerations for the associated nonlinear subproblem permit the calculation of an explicit solution representation.
In Figures~\ref{fig:Figure4} and~\ref{fig:Figure5}, we constrast standard implementations with severe order reductions for a complex Ginzburg--Landau equation~\eqref{eq:CGL}, a Gross--Pitaevskii equation~\eqref{eq:GPE}, and a related parabolic equation~\eqref{eq:GPEParabolic} to a correct implementation, where the nonstiff orders are retained. 
In Figure~\ref{fig:Figure6}, we demonstrate the implications of the stability conditions~\eqref{eq:StabilityConditions} for three-dimensional complex Ginzburg--Landau equations~\eqref{eq:CGL}.
When increasing the imaginary part of a constant, the real splitting methods by Yoshida (Method~3) as well as the schemes involving complex coefficients $(a_j)_{j=1}^s$ (Methods~4 and~8) become unstable, whereas the schemes involving non-negative coefficients $(a_j)_{j=1}^s$ (Methods~10, 11, 12) remain stable. 
Overall, the sixth-order schemes are again superior in efficiency when lower tolerances are desirable. 

\MyParagraph{Long-term integrations}
The formation of quasicrystalline patterns in long-term integrations is illustrated in Figure~\ref{fig:Figure3}.
For visualisations, see also
\begin{center}
{\scriptsize
\url{techmath.uibk.ac.at/mecht/MyHomepage/Research/Movie2024Quasicrystal1.m4v} \\[1mm]
\url{techmath.uibk.ac.at/mecht/MyHomepage/Research/Movie2024Quasicrystal2.m4v} \\
}
\end{center}
Long-term integrations of one- and two-dimensional complex Ginzburg--Landau equations modelling nonlinear waves are illustrated in Figure~\ref{fig:Figure7}, see also 
\begin{center}
{\scriptsize
\url{techmath.uibk.ac.at/mecht/MyHomepage/Research/Movie2024GinzburgLandau1.m4v} \\[1mm]
\url{techmath.uibk.ac.at/mecht/MyHomepage/Research/Movie2024GinzburgLandau2.m4v} \\
}
\end{center}
\section*{Acknowledgements}
Part of this work was developed during a research stay at the Wolfgang Pauli Institute Vienna; the authors are grateful to the director Norbert Mauser and the staff
members for their support and hospitality.
This work has been funded by Ministerio de Ciencia e Innovaci{\'o}n (Spain) through projects PID2022-136585NB-C21 and PID2022-136585NB-C22, MCIN/AEI/10.13039/501100011033/FEDER, UE. 
Sergio Blanes and Fernando Casas acknowledge the support of the Conselleria d'Innovaci{\'o}, Universitats, Ci{\`e}ncia i Societat Digital from the Generalitat Valenciana (Spain) through project CIAICO/2021/180.
\newcommand{\MyBibitem}[2]{\bibitem{#1}}
\newcommand{\MyReference}[3]{\textsc{#1}. \emph{#2}. #3.}

\end{document}